\documentclass[a4paper]{amsartmod1}

\usepackage{amsmath,amssymb,amsfonts}
\usepackage{exscale,url,doi,graphicx,psfrag,xcolor}
\usepackage[linesnumbered,ruled,vlined,algosection]{algorithm2e}

\definecolor{color1}{rgb}{0,0,1}
\hypersetup{pdfborder={1 0 0}}

\def\diag{\ensuremath{\mathrm{diag}}}
\def\la{\ensuremath{\lambda}}
\def\RR{\ensuremath{\mathbb R}}
\def\CC{\ensuremath{\mathbb C}}
\def\spanl{\ensuremath{\mathrm{span}}}

\newtheorem{theorem}{Theorem}[section]
\newtheorem{lemma}[theorem]{Lemma}

\theoremstyle{definition}
\newtheorem{definition}[theorem]{Definition}
\theoremstyle{remark}
\newtheorem{remark}[theorem]{Remark}
\numberwithin{equation}{section}
\numberwithin{table}{section}
\numberwithin{figure}{section}

\begin{document}

\title[Improving preconditioned CG-like eigensolvers]
{Toward genuine efficiency and cluster robustness \\ of preconditioned CG-like eigensolvers}

\author[M.~Zhou]{Ming Zhou}
\address{Universit\"at Rostock, Institut f\"ur Mathematik, 
        Ulmenstra{\ss}e 69, 18055 Rostock, Germany}
\email{ming.zhou at uni-rostock (dot) de}

\author[K.~Neymeyr]{Klaus Neymeyr}
\address{Universit\"at Rostock, Institut f\"ur Mathematik, 
        Ulmenstra{\ss}e 69, 18055 Rostock, Germany}
\email{klaus.neymeyr at uni-rostock (dot) de}


\keywords{Hermitian eigenvalue problems, conjugate gradient,
single-vector eigensolvers, cluster robustness, 
two-term recurrences, trial subspace augmentation. \hfill January 7, 2026.} %

\begin{abstract}
The performance of eigenvalue problem solvers (eigensolvers) depends on various factors
such as preconditioning and eigenvalue distribution. Developing stable and rapidly converging
vectorwise eigensolvers is a crucial step in improving the overall efficiency
of their blockwise implementations. The present paper is concerned with
the locally optimal block preconditioned conjugate gradient (LOBPCG) method
for Hermitian eigenvalue problems, and motivated by two recently proposed
alternatives for its single-vector version LOPCG. A common basis of these eigensolvers
is the well-known CG method for linear systems. However, the optimality
of CG search directions cannot perfectly be transferred to CG-like eigensolvers.
In particular, while computing clustered eigenvalues, LOPCG and its alternatives
suffer from frequent delays, leading to a staircase-shaped convergence behavior
which cannot be explained by the existing estimates.
Keeping this in mind, we construct a class of cluster robust vector iterations
where LOPCG is replaced by asymptotically equivalent two-term recurrences
and the search directions are timely corrected by selecting a far previous iterate
as augmentation. The new approach significantly reduces
the number of required steps and the total computational time.
\end{abstract}

\maketitle
\pagestyle{myheadings}
\thispagestyle{plain}

\section{Introduction}\label{sec1}

The conjugate gradient method (CG) for solving linear systems of Hermitian
positive semidefinite matrices has stimulated the development of efficient
and memory-saving iterative solvers for more general non-quadratic optimization problems.
In the context of computing extreme eigenvalues of large and sparse
Hermitian matrices or matrix pairs, there exist various so-called CG-like eigensolvers.
Two typical ones are the locally optimal block preconditioned conjugate gradient method 
using a CG-based three-term recurrence (LOBPCG)~\cite{KNY2001}
and the Jacobi-Davidson method equipped with CG in its inner iterations
where the stopping criterion depends on Rayleigh quotient residuals (JDCG)~\cite{NY2002}. 
Although not revealed in the title, preconditioning is also an important feature
of JDCG, and these two methods have been compared by numerical tests
in a number of studies leading to somewhat disputable comments~\cite{AHLT2005,STA2007,OVT2008}.
Ever more frequently, their convergence theory and alternative implementations
are considered in further research~\cite{HL2006,NEY2009,DSYG2018,WXS2019}. The present paper is motivated by
the recent contributions~\cite{SCB2025,SC2025} concerning the single-vector version LOPCG
of LOBPCG. Therein two novel methods with similar convergence behavior and
memory requirement are provided together with convergence estimates
in the style of a conjecture from~\cite{KNY2001}. The question arises how far
the convergence analysis and the performance improvement of LOPCG 
can take advantage of these results. In particular, we tend to find out
whether the associated estimates are meaningful in most cases
and whether LOPCG can be substantially accelerated by certain
non-block iterations without essential subspace expansion.
An ideal acceleration by possibly adaptive vector iterations would be, after combined with deflation,
superior to other eigensolvers while computing a multitude of (tightly) clustered eigenvalues.

The LOPCG method was introduced in~\cite{KNY2001} in the form of an algorithm
for solving generalized eigenvalue problems of real symmetric matrices and can easily be applied
to Hermitian matrices. Here we focus on a problem setting that is mostly used
in the cited references, namely,
\begin{equation}\label{evp}
 Ax=\la Mx, \quad A,M\in\CC^{n \times n}\ \mbox{Hermitian}, \quad
 M\ \mbox{positive definite}.
\end{equation}
A more general formulation is introduced in the appendix
concerning partial differential operators, Hermitian definite matrix pencils \cite{KPS2014}
and the linear response eigenvalue problem \cite{BL2013}. The fundamental idea of LOPCG
was even earlier discussed in \cite{KNY1987,KNY1998}. Considering that in \eqref{evp}
the smallest eigenvalue $\la_1$ and an associated eigenvector $x_1$ are of interest,
one can heuristically begin with the ideal case that $\la_1$ is available
(if $M$ is an identity matrix, $\la_1$ could be, but not necessarily efficiently,
approximated by a classical implementation of the Lanczos method
without determining approximate eigenvectors~\cite{CW2002}).
Thus it is possible to compute $x_1$ by solving the homogeneous linear system
\begin{equation}\label{lse}
 A_{\la_1}x_1=0 \quad\mbox{for}\quad A_{\la_1}=A-\la_1M
\end{equation}
with CG thanks to the positive semidefiniteness of $A_{\la_1}$.
The well-known CG form with two-term recurrences of iterates and search directions
can be transformed into a three-term recurrence so that the new iterate $x^{(i+1)}$
depends on its predecessors $x^{(i-1)}$, $x^{(i)}$ and the residual $A_{\la_1}x^{(i)}$.
Returning to the practical case that $\la_1$ is unknown, one has to use an approximate eigenvalue
instead and determine the recurrence parameters in another way, e.g., implicitly
by the Rayleigh-Ritz procedure as suggested in LOPCG. It is notable that
for certain matrices this realistic approach converges almost as fast as the heuristic one,
even if the utilized approximate eigenvalues are not close to $\la_1$ in the first steps~\cite{KNY2001}.
In this sense, LOPCG is able to achieve a global optimality among eigensolvers with
the same initial setting and memory requirement. However, its block version
LOBPCG was not always the clear winner in numerical comparisons, especially for
a large number of target eigenvalues with possibly small gaps~\cite{AHLT2005,STA2007,OVT2008}.
Sometimes iterations based on two-term recurrences of CG can provide slightly
better performance, and restarted generalized Davidson or Jacobi-Davidson schemes prevail at least
in the case of less accurate preconditioners. Of course, this occasional suboptimality of LOBPCG
is more or less problem-dependent and does not contradict its popularity and flexibility;
cf.~the applications~\cite{BL2013,LSE2013,KPS2014,LT2015,AET2017,SET2018}.
Recently some numerical tests directly concerning
the single-vector form LOPCG are presented in~\cite{SCB2025,SC2025}
which point out problematic cases and inspire us to accelerate LOPCG economically, e.g.,
with at most four-dimensional trial subspaces.

For this purpose, we need to classify the test problems with respect to the eigenvalue distribution
and the preconditioning quality on the basis of the available convergence theory of LOPCG.
The first-ever estimate is indeed a conjecture from~\cite{KNY2001} via
the observation that expanding the trial subspace of LOPCG backwards by several more
previous iterates cannot significantly accelerate the convergence. In such cases, LOPCG can be
considered as an almost optimally restarted generalized Davidson iteration and, in particular,
analyzed in terms of an Invert-Lanczos process when preconditioning is extremely fine~\cite{NEY2009}.
The resulting bounds indicate that the convergence rate is mainly influenced by
the gap $\la_2-\la_1$ between the two smallest distinct eigenvalues. We note that the phenomenon
observed in~\cite{KNY2001} is related to the harmless gap in the test problems therein.
In contrast, if $\la_1$ and $\la_2$ are clustered, the eigenvalue approximation by LOPCG
could have a staircase-shaped convergence behavior (cf.~Figure~\ref{fig2} below)
whereas Invert-Lanczos remains robust, and it would be a similar case
for LOBPCG if the cluster size of target eigenvalues exceeds the block size.
Furthermore, for LOPCG with less accurate preconditioners, one has to integrate
certain quality parameters into the analysis. However, this does not mean
a straightforward modification of the estimates for Invert-Lanczos. It is still challenging
to fully derive the convergence factor suggested by the conjecture from~\cite{KNY2001}.
Existing estimates are mostly derived by comparative analyses via simplified iterations
such as the preconditioned steepest descent iteration (PSD) where the iterate $x^{(i-1)}$
is removed from the trial subspace of LOPCG. Among them are the sharp estimates
for PSD with simple step sizes~\cite{KNN2003} or Rayleigh-Ritz step sizes \cite{NEY2012}
where the convergence factors only depend on two or three eigenvalues and
a quality parameter of preconditioning. A similar estimate from \cite{OVT2006s}
additionally uses the current approximate eigenvalue to construct a dynamic convergence factor.
Although these results are not sharp while applied to LOPCG,
they still clearly reflect the fact that an outstanding convergence speed can be guaranteed by
preconditioners well fulfilling the associated quality conditions. Summarizing the above,
challenging cases for LOPCG can principally be built by selecting problems with clustered eigenvalues
and less accurate preconditioners.

To improve the performance of LOPCG in problematic cases, we first construct
a direct modification named LOPCGa where the trial subspace is augmented with
an auxiliary vector which needs to be updated if its angle deviation from the current iterate
is above a threshold. This technique dates back to
an implementation of PSD with implicit deflation~\cite{CBPS2018,ZBCN2023}
and has recently been applied to an eigensolver based on preconditioning and
implicit convexity (EPIC)~\cite{SCB2025}, but not utilized in a related scheme
named Riemannian acceleration with preconditioning (RAP)~\cite{SC2025}.
The dimension of trial subspace reads $3$ in RAP, LOPCG and $4$ in EPIC, LOPCGa.
According to the accompanying and our own numerical tests, EPIC is able to accelerate LOPCG
with respect to the number of required steps and also the total computational time,
whereas RAP is less competitive but still improves PSD, and LOPCGa is mostly more efficient
than EPIC. Nevertheless, EPIC and RAP are parameter-dependent methods and can potentially
be strengthened by advanced strategies in their further development. For LOPCGa, we subsequently
enhance the simple update criterion of the auxiliary vector by detecting peaks
in the convergence history of the residual norm, leading to even better performance.
This encourages us to similarly modify the CG eigensolvers with two-term recurrences
presented in~\cite{OVT2008a} which are asymptotically equivalent to LOPCG.

We structure the remainder of this paper as follows: Section~\ref{sec2} introduces
basic settings and reviews LOPCG concerning its motivation and limitation
in comparison to eigensolvers with similar memory requirement
including the modified method LOPCGa.
Section~\ref{sec3} discusses the state-of-the-art estimates for LOPCG and related eigensolvers.
Section~\ref{sec4} is devoted to appropriate update strategy of the auxiliary vector in LOPCGa
and more efficient schemes with augmented two-term recurrences.
Numerical tests are distributed among these sections.

\section{LOPCG and comparable schemes}\label{sec2}

We consider the generalized eigenvalue problem \eqref{evp} with eigenvalues
$\la_1<\la_2\le\cdots\le\la_n$ and focus on the numerical computation of
an eigenvector $x_1$ associated with $\la_1$ by preconditioned CG-like iterations. The case
$\la_1=\cdots=\la_s<\la_{s+1}$ concerning multiple eigenvalues is neglected
for better readability (a reformulation is straightforward by using projections onto eigenspaces).
The vector iterates approximating $x_1$ are denoted by $x^{(i)}$ with the iteration index $i$.
The corresponding approximate eigenvalues are determined by evaluating the Rayleigh quotient
\begin{equation}\label{rq}
 \la:\CC^n{\setminus}\{0\}\to\RR, \quad \la(x)=\frac{x^*Ax}{x^*Mx}
\end{equation}
where the asterisk ${}^*$ stands for the complex conjugate transpose. The value $\la(x^{(i)})$
is simply denoted by $\la^{(i)}$. In general, any Hermitian positive definite matrix
$T\in\CC^{n \times n}$ can be used as a preconditioner and interpreted as
\begin{equation}\label{prec}
 T \approx A_{\sigma}^{-1} 
 \quad\mbox{with}\quad \kappa=\mbox{cond}_2(T^{1/2}A_{\sigma}T^{1/2})
 \quad\mbox{for}\quad
 A_{\sigma}=A-\sigma M, \quad \sigma<\la_1
\end{equation}
where $\mbox{cond}_2$ denotes the spectral condition number, i.e.,
$\kappa$ coincides with $\alpha_n/\alpha_1$ for the extreme eigenvalues
$\alpha_1$ and $\alpha_n$ of $T^{1/2}A_{\sigma}T^{1/2}$.
The practical construction of $T$ can utilize incomplete factorizations
or approximate solutions of linear systems concerning $A_{\sigma}$. 
The so-called non-preconditioned case means that $T$ is given by the $n \times n$
identity matrix $I$. We note that the non-preconditioned case
can make the analysis easier only for standard eigenvalue problems (i.e. $M=I$),
but not for generalized eigenvalue problems.

\subsection{Benchmark method}\label{sec2a}

The first scheme in our review is a heuristic benchmark suggested by Knyazev
in~\cite[Section~3]{KNY2001} where $\la_1$ is assumed to be available and
the homogeneous linear system \eqref{lse} is to be solved by the standard
preconditioned CG method, e.g., in the form of Algorithm~\ref{apcg}.

\begin{algorithm}
\caption{Heuristic PCG concerning the linear system \eqref{lse}}\label{apcg}
 generate $x^{(0)}$; set $r^{(0)}=-A_{\la_1}x^{(0)}$;
 $w=Tr^{(0)}$; $\gamma^{(0)}=w^*r^{(0)}$; $p^{(0)}=w$\;
\For{$i=0,1,\ldots$ {\rm until convergence}}
{$w=A_{\la_1}p^{(i)}$; \quad $\delta=\gamma^{(i)}/(w^*p^{(i)})$\;
 $x^{(i+1)}=x^{(i)}+\delta p^{(i)}$\;
 $r^{(i+1)}=r^{(i)}-\delta w$\;
 $w=Tr^{(i+1)}$; \quad $\gamma^{(i+1)}=w^*r^{(i+1)}$\;
 $p^{(i+1)}=w+(\gamma^{(i+1)}/\gamma^{(i)})p^{(i)}$\;}
\end{algorithm}

Some convergence properties have been described in~\cite{KNY2001},
however, without an explicit estimate for approximate eigenvalues that can directly
be compared with available estimates for related schemes. We fill this gap
by adding some proofs and remarks. Therein the trivial (and practically rare) case
that $x^{(0)}$ is an eigenvector associated with $\la_1$ or $M$-orthogonal
to the corresponding eigenspace is neglected.

\begin{lemma}\label{lmpcg}
Consider Algorithm~\ref{apcg} with $T$ from \eqref{prec}
and define $B=T^{1/2}A_{\la_1}T^{1/2}$. Let $V$ be an arbitrary
orthonormal basis matrix of the image $\mbox{im}(B)$ of $B$, then the vectors
$\tilde{x}^{(i)}=V^*T^{-1/2}x^{(i)}$, $i=0,1,\ldots$ are CG-iterates for solving
the linear system $\tilde{B}\tilde{x}_1=0$ for $\tilde{B}=V^*BV$.
\end{lemma}
\begin{proof}
Noting that $VV^*$ is the projection matrix of $\mbox{im}(B)$, we have
$B=(VV^*)B(VV^*)$. Then the residual $r^{(i)}$ of $x^{(i)}$
can be converted to the residual $\tilde{r}^{(i)}$ of $\tilde{x}^{(i)}$ by
\[\begin{split}
 & r^{(i)}=-A_{\la_1}x^{(i)}\\
 \Rightarrow\quad
 &T^{1/2}r^{(i)}=-T^{1/2}A_{\la_1}T^{1/2}T^{-1/2}x^{(i)}
  =-BT^{-1/2}x^{(i)}=-VV^*BVV^*T^{-1/2}x^{(i)}\\
 \Rightarrow\quad
 &V^*T^{1/2}r^{(i)}=-V^*VV^*BVV^*T^{-1/2}x^{(i)}=-\tilde{B}\tilde{x}^{(i)}=\tilde{r}^{(i)}.
\end{split}\]
Next, we check the conversions of the terms $\gamma^{(i)}$ and $Tr^{(i)}$
concerning the update of the search direction $p^{(i)}$, i.e.,
\[\begin{split}
 \gamma^{(i)}=r^{(i)}{}^*Tr^{(i)}
 &= (T^{1/2}r^{(i)})^*(T^{1/2}r^{(i)})
   =(-VV^*BVV^*T^{-1/2}x^{(i)})^*(-VV^*BVV^*T^{-1/2}x^{(i)})\\
 &= (-V^*BVV^*T^{-1/2}x^{(i)})^*(-V^*BVV^*T^{-1/2}x^{(i)})
   =\tilde{r}^{(i)}{}^*\tilde{r}^{(i)},\\
 (V^*T^{-1/2})Tr^{(i)} &= V^*T^{1/2}r^{(i)}=\tilde{r}^{(i)}.
\end{split}\]
Thus Algorithm~\ref{apcg} corresponds to a (non-preconditioned) CG-scheme
for solving $\tilde{B}\tilde{x}_1=0$.
\end{proof}

Lemma~\ref{lmpcg} allows applying the well-known convergence estimates
for CG since $\tilde{B}$ is evidently positive definite. Based on~\cite[Theorem~6.29]{SY2003},
it holds that
\begin{equation}\label{lmpcge}
 \frac{\|\tilde{x}^{(i)}\|_{\tilde{B}}}{\|\tilde{x}^{(0)}\|_{\tilde{B}}}
 =\frac{\|0-\tilde{x}^{(i)}\|_{\tilde{B}}}{\|0-\tilde{x}^{(0)}\|_{\tilde{B}}}
 \le \big(\mathcal{C}_i(\tilde{\varphi})\big)^{-1} \quad\mbox{for}\quad
 \tilde{\varphi}=\frac{\tilde{\beta}_m+\tilde{\beta}_1}{\tilde{\beta}_m-\tilde{\beta}_1}
\end{equation}
with the Chebyshev polynomial $\mathcal{C}_i$ of degree $i$ of the first kind
and the smallest eigenvalue $\tilde{\beta}_1$ and
the largest eigenvalue $\tilde{\beta}_m$ of $\tilde{B}$.
Our next task is to extend \eqref{lmpcge} to an explicit estimate for Algorithm~\ref{apcg}.

\begin{theorem}\label{thmpcg}
The approximate eigenvalues $\la^{(i)}=\la(x^{(i)})$, $i=0,1,\ldots$ by Algorithm~\ref{apcg} fulfill
\begin{equation}\label{thmpcge}
 \la^{(i)}-\la_1 \le\,
 \big(\mathcal{C}_i(\varphi)\big)^{-2}
 \frac{\|x^{(0)}\|_M^2}{\|x^{(i)}\|_M^2}\,(\la^{(0)}-\la_1) \quad\mbox{for}\quad
 \varphi=\frac{\eta+1}{\eta-1}, \quad
 \eta=\kappa\,\frac{(\la_n-\la_1)(\la_2-\sigma)}{(\la_2-\la_1)(\la_n-\sigma)}
\end{equation}
where $\kappa$ is the condition number defined in \eqref{prec}.
\end{theorem}
\begin{proof}
The matrix $B=T^{1/2}A_{\la_1}T^{1/2}$ given in Lemma~\ref{lmpcg} is positive semidefinite
and has rank $n-1$ by checking the properties of $A_{\la_1}$ and $T$. Denoting its eigenvalues by
$\beta_1\le\beta_2\le\cdots\le\beta_n$, it holds that $\beta_1=0$, $\beta_j>0$ for $j=2,\ldots,n$.
Moreover, since $\tilde{B}=V^*BV$ is the restriction of $B$ to $\mbox{im}(B)$,
the extreme eigenvalues $\tilde{\beta}_1$ and $\tilde{\beta}_m$ of $\tilde{B}$
coincide with $\beta_2$ and $\beta_n$, respectively.

On the other hand, the eigenvalues of $B=T^{1/2}A_{\la_1}T^{1/2}$ are also those of
the matrix pair $(A_{\la_1},T^{-1})$ and thus related to the Rayleigh quotient
\[\beta(y)=\frac{y^*A_{\la_1}y}{y^*T^{-1}y}\]
which can be rewritten as the product of two other Rayleigh quotients: 
\[\beta(y)=\psi(y)\,\alpha(y)\quad\mbox{for}\quad
 \psi(y)=\frac{y^*A_{\la_1}y}{y^*A_{\sigma}y}
 \quad\mbox{and}\quad \alpha(y)=\frac{y^*A_{\sigma}y}{y^*T^{-1}y}.\]
Correspondingly, we denote the eigenvalues of the matrix pairs
$(A_{\la_1},A_{\sigma})$ and $(A_{\sigma},T^{-1})$ by
$\psi_1\le\psi_2\le\cdots\le\psi_n$ and $\alpha_1\le\alpha_2\le\cdots\le\alpha_n$
(which are all non-negative). Then $\alpha_1 \psi(y) \le \beta(y) \le \alpha_n \psi(y)$ holds,
and the Courant-Fischer principles imply
\[\begin{split}
 \beta_2&=\min_{\dim\mathcal{Y}=2}\max_{y\in\mathcal{Y}{\setminus}\{0\}}\beta(y)
 \ge \alpha_1\min_{\dim\mathcal{Y}=2}\max_{y\in\mathcal{Y}{\setminus}\{0\}}\psi(y) = \alpha_1\psi_2,\\
 \beta_n&=\max_{y\in\CC^n{\setminus}\{0\}}\beta(y)
 \le \alpha_n\max_{y\in\CC^n{\setminus}\{0\}}\psi(y) = \alpha_n\psi_n.
\end{split}\]
Since $\alpha_1$ and $\alpha_n$ are also the extreme eigenvalues of $T^{1/2}A_{\sigma}T^{1/2}$,
it holds that $\kappa=\alpha_n/\alpha_1$ by~\eqref{prec}.
The eigenvalues of $(A_{\la_1},A_{\sigma})$
are those of $A_{\sigma}^{-1}A_{\la_1}$ and can be represented in terms of the eigenvalues
of $(A,M)$. By using the diagonal matrix
$\Lambda=\diag(\la_1,\ldots,\la_n)$ and a matrix $X$ consisting of
associated $M$-orthonormal eigenvectors as columns, we get the factorization
\[A_{\sigma}^{-1}A_{\la_1}
 =\big(X (\Lambda -\sigma I)^{-1}X^*\big)\big(MX (\Lambda -\la_1 I)X^*M\big)
 =X(\Lambda -\sigma I)^{-1}(\Lambda -\la_1 I)X^{-1}\]
so that the eigenvalues of $(A_{\la_1},A_{\sigma})$
are given by the diagonal of $(\Lambda -\sigma I)^{-1}(\Lambda -\la_1 I)$, i.e.,
$\psi_j=(\la_j-\la_1)/(\la_j-\sigma)$ for $j=1,\ldots,n$. Summarizing the above,
\begin{equation}\label{thmpcgf}
 \frac{\tilde{\beta}_m}{\tilde{\beta}_1}=\frac{\beta_n}{\beta_2}
 \le\frac{\alpha_n\psi_n}{\alpha_1\psi_2}
 =\kappa\,\frac{(\la_n-\la_1)(\la_2-\sigma)}{(\la_2-\la_1)(\la_n-\sigma)}=\eta,
\end{equation}
and the Chebyshev term in \eqref{lmpcge} is bounded in terms of
\begin{equation}\label{thmpcge1}
 \tilde{\varphi}=\frac{\tilde{\beta}_m+\tilde{\beta}_1}{\tilde{\beta}_m-\tilde{\beta}_1}
 =\frac{\tilde{\beta}_m/\tilde{\beta}_1+1}{\tilde{\beta}_m/\tilde{\beta}_1-1}
 \ge\frac{\eta+1}{\eta-1}=\varphi>1 \quad\Rightarrow\quad
 \big(\mathcal{C}_i(\tilde{\varphi})\big)^{-1}\le\big(\mathcal{C}_i(\varphi)\big)^{-1}.
\end{equation}

Next, the norm transformation (using again $B=(VV^*)B(VV^*)$)
\[\begin{split}
 \|\tilde{x}^{(i)}\|_{\tilde{B}}^2 & =\tilde{x}^{(i)}{}^*\tilde{B}\tilde{x}^{(i)}
 =(V^*T^{-1/2}x^{(i)})^*(V^*BV)(V^*T^{-1/2}x^{(i)}) \\
 & =x^{(i)}{}^*T^{-1/2}BT^{-1/2}x^{(i)}
 =x^{(i)}{}^*T^{-1/2}(T^{1/2}A_{\la_1}T^{1/2})T^{-1/2}x^{(i)}\\
 &=x^{(i)}{}^*A_{\la_1}x^{(i)}
 =x^{(i)}{}^*Ax^{(i)}-x^{(i)}{}^*\la_1Mx^{(i)}
 =(\la^{(i)}-\la_1)\|x^{(i)}\|_M^2
\end{split}\]
and its counterpart for $x^{(0)}$ complete the proof by combining $\eqref{thmpcge1}$
with \eqref{lmpcge}.
\end{proof}

Some geometric relations between $x^{(i)}$ and the target eigenspace are useful for
extending Theorem~\ref{thmpcg}.

\begin{lemma}\label{lmpcg1}
Let $\mathcal{X}_1$ be the eigenspace associated with $\la_1$. 
Then the iterates $x^{(i)}$, $i=0,1,\ldots$ by Algorithm~\ref{apcg} fulfill
\begin{equation}\label{lmpcg1e}
 x^{(i)}-x^{(0)}\perp_{T^{-1}} \mathcal{X}_1 \quad\mbox{and}\quad
 \|x^{(i)}\|_M^2\ge\mu_1\|x^{(0)}\|_{T^{-1}}^2\big(\cos\angle_{T^{-1}}(x^{(0)},\mathcal{X}_1)\big)^2
\end{equation}
where $\mu_1$ is the smallest eigenvalue of the matrix pair $(M,T^{-1})$
and $\angle_{T^{-1}}$ denotes angles with respect to the inner product induced by $T^{-1}$.
\end{lemma}
\begin{proof}
The CG property
\[\begin{split}
 x^{(i)} & \in x^{(0)}+\spanl\{Tr^{(0)},\,(TA_{\la_1})\,Tr^{(0)},\ldots,(TA_{\la_1})^{i-1}Tr^{(0)}\}\\
 & = x^{(0)}+(TA_{\la_1})\,\spanl\{x^{(0)},\,(TA_{\la_1})\,x^{(0)},\ldots,(TA_{\la_1})^{i-1}x^{(0)}\}
\end{split}\]
ensures $x^{(i)}-x^{(0)}=TA_{\la_1}w$ for certain $w$. Thus it holds for
any $x_1\in\mathcal{X}_1{\setminus}\{0\}$ that
\[x_1^*T^{-1}(x^{(i)}-x^{(0)})=x_1^*T^{-1}TA_{\la_1}w=(A_{\la_1}x_1)^*w=0\]
which yields the orthogonality in \eqref{lmpcg1e}. We further denote by $x_1$ the $T_{-1}$-projection
of $x^{(0)}$ onto $\mathcal{X}_1$, then this orthogonality implies
\[\begin{split}\cos\angle_{T^{-1}}(x^{(0)},\mathcal{X}_1)
 &=\cos\angle_{T^{-1}}(x^{(0)},x_1)
   =\frac{x_1^*T^{-1}x^{(0)}}{\|x_1\|_{T^{-1}}\|x^{(0)}\|_{T^{-1}}}\\
 &=\frac{x_1^*T^{-1}x^{(i)}}{\|x_1\|_{T^{-1}}\|x^{(0)}\|_{T^{-1}}}
   \le\frac{\|x_1\|_{T^{-1}}\|x^{(i)}\|_{T^{-1}}}{\|x_1\|_{T^{-1}}\|x^{(0)}\|_{T^{-1}}}
   =\frac{\|x^{(i)}\|_{T^{-1}}}{\|x^{(0)}\|_{T^{-1}}}.
\end{split}\]
Moreover, by noting that $\|x^{(i)}\|_M^2/\|x^{(i)}\|_{T^{-1}}^2$ is the Rayleigh quotient
of the matrix pair $(M,T^{-1})$ evaluated at $x^{(i)}$, a lower bound is given by $\mu_1$.
Thus \\[1ex] \mbox{\quad}\hspace{3cm}
$\displaystyle \|x^{(i)}\|_{M}^2 \ge \mu_1\|x^{(i)}\|_{T^{-1}}^2
 \ge \mu_1\|x^{(0)}\|_{T^{-1}}^2\big(\cos\angle_{T^{-1}}(x^{(0)},\mathcal{X}_1)\big)^2.$
\end{proof}

The $T^{-1}$-orthogonality in \eqref{lmpcg1e} indicates that the $T^{-1}$-angle between
$x^{(i)}$ and $\mathcal{X}_1$ can measure the convergence of Algorithm~\ref{apcg},
whereas the inequality in \eqref{lmpcg1e} can directly be added to \eqref{thmpcge} to eliminate
the term $\|x^{(i)}\|_M^2$.

\begin{remark}
In the special case $T=M=I$, combining \eqref{lmpcg1e} with \eqref{thmpcge}
implies an estimate for the Lanczos method.
Therein $\mu_1=1$, and \eqref{lmpcg1e} gives
\begin{equation}\label{lmpcg1e1}
 \|x^{(i)}\|_2^2\ge\|x^{(0)}\|_2^2\big(\cos\angle_2(x^{(0)},\mathcal{X}_1)\big)^2
\end{equation}
(where the subscript ${}_2$ denotes Euclidean terms). Moreover, going through the proof
of Theorem~\ref{thmpcg} with $T=I$ leads to $\eta=(\la_n-\la_1)/(\la_2-\la_1)$. Then
\eqref{thmpcge} becomes, by using \eqref{lmpcg1e1},
\begin{equation}\label{lmpcg1e2}
 \la^{(i)}-\la_1 \le\,
 \big(\mathcal{C}_i(\varphi)\big)^{-2}
 \big(\cos\angle_2(x^{(0)},\mathcal{X}_1)\big)^{-2}(\la^{(0)}-\la_1) \quad\mbox{for}\quad
 \varphi=\frac{\eta+1}{\eta-1}=1+2\,\frac{\la_2-\la_1}{\la_n-\la_2}
\end{equation}
which differs from its traditional counterpart from \cite[Theorem~2]{SY1980}
(by changing the eigenvalue arrangement) only in terms depending on
the initial guess. Subsequently, \eqref{lmpcg1e2} is directly applicable to
the Lanczos method since for $T=M=I$ Algorithm~\ref{apcg}
generates the same series of Krylov subspaces with respect to $A$,
and $\la^{(i)}$ is not smaller than its counterpart in the Lanczos method,
i.e., the smallest Ritz value associated with the current Krylov subspace.
In addition, \eqref{lmpcg1e2} gives a tighter bound than the traditional one
due to $\la^{(0)}-\la_1 \le (\la_n-\la_1)\big(\sin\angle_2(x^{(0)},\mathcal{X}_1)\big)^2$.
Analogously, in the special case $T=A_{\sigma}^{-1}$, Algorithm~\ref{apcg} produces Krylov subspaces
with respect to $A_{\sigma}^{-1}M$ so that \eqref{thmpcge} implies an estimate
for the Invert-Lanczos method which is comparable with the one from~\cite{NEY2009}.
The case $T=I$, $M \neq I$ is, however, not related to the Lanczos and  Invert-Lanczos methods.
\end{remark}

\begin{remark}
In the final phase where $\la^{(i)}$ is close to $\la_1$,
Algorithm~\ref{apcg} behaves like the generalized Davidson method
so that a comparative analysis is conceivable.
However, combining \eqref{lmpcg1e} with \eqref{thmpcge}
could cause a substantial overestimation for less accurate preconditioners where
the sharpness depends on the choice of $\sigma<\la_1$ in \eqref{prec}.
Although the pseudoinverse of $A_{\la_1}$ leads to a one-step convergence,
selecting $\sigma$ from the near vicinity of $\la_1$ is not reasonable for constructing
practical preconditioners. For instance, as far as $A$ is known to be positive definite,
one usually sets $\sigma=0$ so that $A_\sigma=A$ and
$\kappa=\mbox{cond}_2(T^{1/2}AT^{1/2})$; cf.~\cite{KNY2001,SCB2025}.
Moreover, the multi-step convergence factor $\big(\mathcal{C}_i(\varphi)\big)^{-2}$ for
$\varphi=(\eta+1)/(\eta-1)$ can be interpreted by an asymptotic average
single-step convergence factor of the form $\big((\sqrt{\eta}-1)/(\sqrt{\eta}+1)\big)^2$
based on equivalent reformulations of Chebyshev polynomials. This interpretation
is problematic for the first steps by noting that the first CG-step is still a simple gradient iteration
so that the convergence factor therein should be close to $\big((\eta-1)/(\eta+1)\big)^2$.
We continue this discussion in Section~\ref{sec3}.
\end{remark}

\subsection{Practical eigensolvers}\label{sec2b}

The only term in Algorithm~\ref{apcg} that hinders a practical implementation
is $A_{\la_1}$ since the eigenvalue $\la_1$ is actually the one to be computed.
A simple way out is replacing $\la_1$ by a readily obtainable approximation,
e.g., the Rayleigh quotient value $\la^{(i)}=\la(x^{(i)})$. Of course,
this replacement alone would deteriorate the convergence, since the residual
and the step size need to be redefined. By considering that $\la_1$ is the minimum
of the Rayleigh quotient \eqref{rq}, a scaled gradient can be used as the residual.
Moreover, the global minimization of \eqref{rq} can be achieved by computing
Ritz values for a series of updated subspaces. 

\begin{definition}
We denote by
\begin{equation}\label{rrw}
x' \ \ \xleftarrow{\mathrm{RRw}} \ \ x+\,\mathcal{U} \ \ : \ \
x' \in  x+\,\mathcal{U} \ \ \mbox{such that} \ \
\la(x')=\min_{w\in\spanl\{x,\,\mathcal{U}\}\setminus\{0\}}\la(w)
\end{equation}
the Rayleigh-Ritz procedure aiming at the smallest Ritz value
but with a weighted output with respect to the first basis vector.
\end{definition}

By using \eqref{rrw}, one can modify Lines 3 to 5 in Algorithm~\ref{apcg} as
\[x^{(i+1)} \ \ \xleftarrow{\mathrm{RRw}} \ \ x^{(i)}+\spanl\{p^{(i)}\};
 \quad r^{(i+1)}=Ax^{(i+1)}-\la^{(i+1)}Mx^{(i+1)};\]
to maintain the convergence toward $\la_1$. In addition,
some other formulas for the search directions $p^{(i)}$
result in better performance \cite{OVT2008}. This kind of modification (or the block form) is, however,
not as popular as the LO(B)PCG method. We discuss its possible improvement in Section~\ref{sec4}.

A characteristic feature of LOPCG is to determine all step sizes implicitly within the Rayleigh-Ritz procedure.
In Algorithm~\ref{apcg}, by combining Lines 6, 7 with Line 4 in the next step,
$x^{(i+2)}$ is evidently a linear combination of $x^{(i+1)}$, $Tr^{(i+1)}$ and $p^{(i)}$.
This suggests a three-dimensional trial subspace and leads to
a basic version of LOPCG as shown in Algorithm~\ref{alopcg}. Therein the iteration index
of $p^{(i)}$ is shifted, and the step sizes in Algorithm~\ref{apcg}
are replaced by certain by-products for determining $x^{(i+1)}$ as a Ritz vector.
The trial subspace in Algorithm~\ref{alopcg} is mostly more stable
than its theoretically equivalent form using $x^{(i-1)}$ instead of $p^{(i)}$.
Besides this well-known fact mentioned in \cite{KNY2001}, it seems to be less known
that the weighted forms $x^{(i+1)}=x^{(i)}+\cdots$ (consistent with the heuristic PCG)
and $x^{(i+1)}=Tr^{(i)}+\cdots$ (suggested in \cite{KNY2001})
without any orthogonalization could show clearly different performance.
Indeed, the practical code \cite[\texttt{lobpcg.m}]{KALO2007} uses a partial orthogonalization
where $Tr^{(i)}$ is orthogonalized against $x^{(i)}$. Moreover, $p^{(i)}$ will be removed
from the trial space if the underlying Gram matrix is ill-conditioned.
Another suggestion for stabilization~\cite{HL2006,DSYG2018}
is to apply a full $M$-orthonormalization of the three basis vectors, yet this strategy
increases the computational time per step and occasionally,
according to our recent tests, even the number of required steps; see Section~\ref{sec2c}.
The two-term CG-schemes from~\cite{OVT2008} seem not to suffer from such instability.

\begin{algorithm}
\caption{LOPCG for computing $(\la_1,x_1)$ in \eqref{evp}}\label{alopcg}
 generate $x^{(0)}$; set $p^{(0)}=0$\;
\For{$i=0,1,\ldots$ {\rm until convergence}}
{$\la^{(i)}=\la(x^{(i)})$; \quad $r^{(i)}=Ax^{(i)}-\la^{(i)}Mx^{(i)}$\;
 $x^{(i+1)} \ \ \xleftarrow{\mathrm{RRw}} \ \ x^{(i)}+\spanl\{Tr^{(i)},\,p^{(i)}\}$ as in \eqref{rrw}\;
 $p^{(i+1)}=x^{(i+1)} -  x^{(i)}$\;}
\end{algorithm}

The performance of LOPCG has been shown in \cite[Section~7]{KNY2001}
with diagonal $A$ and $M=I$. The preconditioner $T$
is constructed by random diagonal matrices multiplied by orthogonal matrices.
The smallest two eigenvalues are well separated: $\la_1=1$ and $\la_2=2$.
The conclusion is that LOPCG performs similarly to the heuristic PCG
for both ideal $T$ with $\kappa=4$ and less accurate $T$ with $\kappa=1000$
with respect to the setting \eqref{prec} for $\sigma=0$.
It is remarkable that although the first approximate eigenvalues $\la^{(i)}$
are usually far from $\la_1$ at the beginning, the goal interval $(\la_1,\la_2)$
can still be reached within a few steps. This also gives a reason why
some results in the available convergence theory simply begin with the assumption
$\la^{(0)}<\theta$ for certain $\theta\in(\la_1,\la_2)$.
Another fact is that LOPCG significantly accelerates
the preconditioned steepest descent iteration (PSD) whose trial subspace
is simply $\spanl\{x^{(i)},\,Tr^{(i)}\}$, analogously to
the comparison between their counterparts for solving linear systems.

In a number of performance studies, LOPCG and its block extension LOBPCG are numerically
compared with other practical eigensolvers, however, not always with similar memory requirement. 
For instance, in \cite{STA2007,WXS2019}, LOPCG and two-block LOBPCG
are considered together with Davidson-type methods which use up to 18-dimensional
trial subspaces, even for computing one extreme eigenpair.
We note that in the situations where LOPCG converges dramatically slower than
Davidson-type methods, its convergence behavior also significantly differs from
that of the heuristic PCG. In our opinion, a fairly meaningful comparative analysis
of LOPCG should focus on eigensolvers using trial subspaces with dimensions 2 to 4,
e.g., the two-term CG-schemes from \cite{OVT2008} mentioned above
as well as the recently developed methods EPIC \cite{SCB2025} and RAP \cite{SC2025}.

We briefly introduce some features of EPIC and RAP since they are less related to the heuristic PCG,
but rather derived from minimizing certain functions other than the Rayleigh quotient
associated with the considered eigenvalue problem.
The approach in EPIC (where $T$ is denoted by $T^{-1}$)
utilizes an auxiliary vector $q$ and considers $M$-normalized vectors $x$
with $q^*Mx>0$. Then an intermediate function defined for the representations of such $x$
underlies the further construction by applying the convex optimization theory.
Apart from differences in detail, the setting for preconditioners is indeed similar to \eqref{prec},
but with respect to the Hessian of the intermediate function. The associated condition number
is represented in the form $L/\mu$ and shown in~\cite[Corollary~2.1]{SCB2025}
to be close to the parameter $\eta$ in \eqref{thmpcge}
for $\sigma=0$ with a distance depending on $\la(q)-\la_1$.
Interestingly, the parameters $\mu$ and $L$ are also involved in
the algorithm of EPIC, and thus need to be estimated before a practical implementation.
Yet the choice $\mu=L=6$ in~\cite{SCB2025} is quite empirical, since this corresponds to
a perfect preconditioner which allows a one-step convergence, but the preconditioners
chosen for the numerical tests therein are evidently less accurate ones.
Nevertheless, according to our additional attempt, using the exact $\eta$ as $L/\mu$
does not lead to a faster convergence than the above choice; see Section~\ref{sec2c}.
On the other hand, the intermediate function in RAP (where $T$ is denoted by $B^{-1}$)
is related to the Rayleigh quotient associated with the matrix pair
$(T^{1/2}AT^{1/2},T^{1/2}MT^{1/2})$ and does not depend on auxiliary vectors.
The condition number of preconditioning also uses the form $L/\mu$
and asymptotically coincides with $\eta$ in \eqref{thmpcge}; cf.~\cite[Corollary~16]{SC2025}.
The parameter choice $\mu=8$, $L=50$ in~\cite{SC2025}
seems to be determined by a concrete numerical example,
but contradicts the condition $L \ge 9\mu$ in~\cite[Proposition~24, Theorem~26]{SC2025}.

In direct comparison to LOPCG, the approaches in EPIC and RAP suggest novel trial subspaces
of similar dimensions (4 and 3). According to the accompanying and our own numerical tests, 
EPIC can be superior to LOPCG whereas RAP mostly converges slower than LOPCG but still accelerates PSD.
This inspires us to construct an efficient modification of LOPCG whose
trial subspace is four-dimensional just like that of EPIC. However, as reported in \cite{KNY2001},
adding directly previous iterates to the trial subspace of LOPCG
does not necessarily lead to a substantial reduction of the number of required steps.
We denote by LOPCGx the simplest version of this type, namely, with
$\spanl\{x^{(i)},\,Tr^{(i)},p^{(i)},x^{(i-2)}\}$ as the trial subspace,
and include LOPCGx in numerical comparisons in Section~\ref{sec2c}.

For the purpose of efficiently modifying LOPCG, we note that
the auxiliary vector $q$ in EPIC is initialized by the initial guess $x^{(0)}$ and that EPIC
is restarted with $q=x^{(i)}$ as soon as $|q^*Mx^{(i)}|$ falls below the empirical value $0.5$.
This restart is necessary since the first $q$ is usually not a good approximate eigenvector
as described in the theoretical construction. Although the criterion was given without explanation,
it is evident that the term $|q^*Mx^{(i)}|$ is the cosine of the $M$-angle
between $x^{(i)}$ and $\spanl\{q\}$ since the EPIC iterates are always $M$-normalized.
Thus the essential information is that an appropriate $q$ should keep a small angle deviation
from the current iterate. The practical effectiveness of $q$ immediately motivates
our first modification of LOPCG called ``augmented LOPCG'' (LOPCGa) where an auxiliary vector $a$
is added to the trial subspace and updated by the criterion $|\cos\angle_M(a,x^{(i)})|<\tau$
with a threshold $\tau=0.7$\,; see Algorithm~\ref{alopcga}.
However, we prefer to interpret $a$ as a far previous iterate since $a$ is, unlike $q$, not used to
generate alternative basis vectors of the trial subspace. In other words, LOPCGa is
still structurally similar to LOPCG and the heuristic PCG.

\begin{algorithm}
\caption{Augmented LOPCG (LOPCGa) for computing $(\la_1,x_1)$ in \eqref{evp}}\label{alopcga}
 generate $x^{(0)}$; set $a=x^{(0)}$\;
\For{$i=0,1,\ldots$ {\rm until convergence}}
{$\la^{(i)}=\la(x^{(i)})$; \quad $r^{(i)}=Ax^{(i)}-\la^{(i)}Mx^{(i)}$; \quad $w=Tr^{(i)}$\;
 \eIf{$i>0$}
 {\eIf{$|\cos\angle_M(a,x^{(i)})|<\tau$}
 {$a=x^{(i)}$; \quad $\mathcal{U}=\spanl\{w,\,p^{(i)}\}$\;}
 {$\mathcal{U}=\spanl\{w,\,p^{(i)},a\}$\;}
  \ $W \ \ \xleftarrow{\mbox{\scriptsize $M$-orthogonalization}} \ \ \spanl\{x^{(i)},\,\mathcal{U}\}$\;
  \If{$W^*MW$ {\rm ill-conditioned}}
  {$\mathcal{U}=\spanl\{w\}$ \ \
   or \ \ $\mathcal{U}=\spanl\{w,\,p^{(i)}\}$\;}}
 {$\mathcal{U}=\spanl\{w\}$\;}
 $x^{(i+1)} \ \ \xleftarrow{\mathrm{RRw}} \ \ x^{(i)}+\,\mathcal{U}$ as in \eqref{rrw};
 \quad $p^{(i+1)}=x^{(i+1)} -  x^{(i)}$\;}
\end{algorithm}

Another important feature of LOPCGa
is inspired by \cite{KALO2007} concerning stabilization via subspace reduction.
As described in Algorithm~\ref{alopcga}, an $M$-orthogonalization in the trial subspace 
decides whether the auxiliary vector $a$ and the search direction $p^{(i)}$
need to be dropped in the current step. More precisely for Lines 10 and 11,
we denote the diagonal entries of  the corresponding (diagonal) Gram matrix $W^*MW$
by $\delta_1,\ldots,\delta_k$ ($k\in\{3,\,4\}$) and utilize a threshold $\gamma=\texttt{1e26}$.
Then the trial subspace is reduced to $\spanl\{x^{(i)},w\}$
in the case $\delta_1>\gamma\delta_3$ or to $\spanl\{x^{(i)},w,\,p^{(i)}\}$
provided that $k=4$ and $\gamma\delta_3\ge\delta_1>\gamma\delta_4$.
This normalization-free technique makes LOPCGa more stable than
its first version using ordinary orthonormalization in the sense of smoother reduction
of the residual norm in the final phase.
Before attempting to construct further modifications, we report the performance of LOPCGa
within a few numerical examples.

\subsection{Numerical comparisons}\label{sec2c}

We select three examples to indicate that LOPCGa can significantly accelerate LOPCG
when there is still space for improving it.
\begin{equation}\label{sp1}
{\small
\begin{array}{cclll}
\mbox{(a)} & A=\texttt{boneS01},\quad M=I, & n=127224,
 & \la_1\approx\texttt{2.847268e-3}, & \la_2\approx\texttt{6.591631e-2}\\[1ex]
\mbox{(b)} & A=\texttt{finan512},\quad M=I, & n=74752,
 & \la_1\approx\texttt{9.474684e-1}, & \la_2\approx\texttt{9.502419e-1}\\[1ex]
\mbox{(c)} & (A,M) \mbox{ by finite elements}, & n=178820,
 & \la_1\approx\texttt{1.973967e1}, & \la_2\approx\texttt{1.973975e1}
\end{array}}
\end{equation}
The examples (a) and (b) listed in \eqref{sp1}
are selected from the SuiteSparse Matrix Collection (University of Florida)
and also used for numerical tests in \cite{SCB2025}
where LOPCG occasionally converges clearly slower than EPIC.
The example (c) is derived from a finite elements discretization
of the Laplacian eigenvalue problem on the rectangle domain
$[0,\,2]{\times}[0,\,1]$ with the slit $\{1\}{\times}[0.1,\,0.9]$
and homogeneous Dirichlet boundary conditions.
Therein the two smallest eigenvalues are tightly clustered.

\begin{figure}[htbp]
\begin{center}
\includegraphics[trim={4.5cm 15cm 4.5cm 4.5cm}, clip, width=\textwidth]{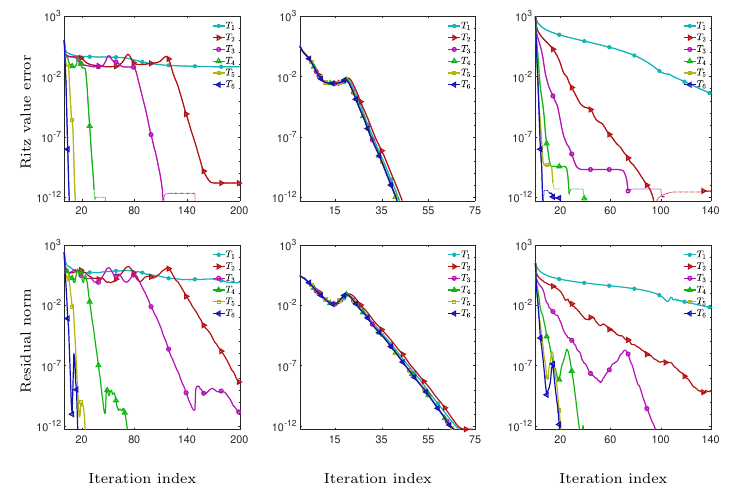}
\end{center}
\par\vskip -2ex
\caption{\small Selecting typical cases by observing
the convergence history of the heuristic PCG (Algorithm~\ref{apcg})
for examples \eqref{sp1} and preconditioners in Table~\ref{tab}.}
\label{figptf}
\end{figure}

For each example, we generate six preconditioners by the function
\texttt{ichol} in Matlab, namely, $T_1$ with \texttt{nofill} and $T_i$ for $i{\,=\,}2,\ldots,6$
with \texttt{ict} using \texttt{droptol}$=$\verb|1e-i| except $T_2$ for (\ref{sp1},\,a) 
which uses \texttt{droptol}$=$\verb|2e-3| due to feasibility.
The associated condition numbers defined in \eqref{prec}
with $\sigma=0$ are listed in Table~\ref{tab}.
In addition, we test the heuristic PCG (Algorithm~\ref{apcg} using $\la_1$ by \texttt{eigs})
for these preconditioners with the $M$-normalized \texttt{ones} vector as initial guess.
The convergence history is illustrated in Figure~\ref{figptf}
with respect to the Ritz value error $\la^{(i)}-\la_1$
and the (relative) residual norm $\|r^{(i)}\|_2/\|x^{(i)}\|_M$.
The performance is clearly related to the conditioner numbers as analyzed
in Section~\ref{sec2a}, apart from the fact that the convergence delay
before the final phase cannot be reflected by any available estimates.
In particular, the residual norm reduction can strongly oscillate for less accurate
preconditioners since the input value of $\la_1$ still differs from the exact eigenvalue.
As a way out, we update this value by the current $\la^{(i)}$ in the final phase,
leading to weaker oscillation as seen for $T_3$ and $T_4$ in the lower left subplot.
This also indicates that a practical eigensolver with unstable final phase can be refined
by switching to the heuristic PCG as soon as an unusual oscillation is detected.
Concerning the Ritz value error, we note that $\la^{(i)}$ can fall below the input value of $\la_1$
so that $\la^{(i)}-\la_1$ cannot be displayed by log scaling. In this case,
we draw $|\la^{(i)}-\la_1|$ by dotted curves; see the upper left/right subplot.
A further fact is that the convergence of $\la^{(i)}$ in the heuristic PCG
is not necessarily monotonic in contrast to iterations using the Rayleigh-Ritz procedure.
Following the above preparation, we select six combinations for
numerical comparisons between practical eigensolvers; see Table~\ref{tab1}.

\begin{table}
\caption{\small Condition numbers of preconditioners.}\label{tab}
\begin{tabular}{c|cccccc}
\hline 
Example & $T_1$ & $T_2$ & $T_3$ & $T_4$ & $T_5$ & $T_6$ \\
\hline
 (\ref{sp1},\,a) 
 & \texttt{4.61e6}  & \texttt{1.15e5}  & \texttt{5.98e4}  & \texttt{4.75e3} & \texttt{3.19e2} & \texttt{6.86e0} \\
 (\ref{sp1},\,b) 
 & \texttt{1.59e0}  & \texttt{1.53e0}  & \texttt{1.07e0}  & \texttt{1.01e0} & \texttt{1.00e0} & \texttt{1.00e0} \\
 (\ref{sp1},\,c) 
 & \texttt{9.31e3}  & \texttt{4.75e2}  & \texttt{5.78e1}  & \texttt{7.04e0} & \texttt{1.45e0} & \texttt{1.02e0} \\
\hline
\end{tabular}
\end{table}

\begin{table}
\caption{\small Parameters for numerical comparisons.}\label{tab1}
\begin{tabular}{c|cccccc|l}
\hline 
Example, $T_i$ & droptol & $\kappa$ & $\eta$ & $\beta_{\min}$ & $\beta_{\max}$ & seed & Figure \\
\hline
 (\ref{sp1},\,a), $T_4$ & \texttt{1e-4}
 & \texttt{4.75e3}  & \texttt{4.96e3}  & \texttt{2.31e-4} & \texttt{1.15e0} & 20 & \,\ref{fig1} left \\
 (\ref{sp1},\,a), $T_3$ & \texttt{1e-3}
 & \texttt{5.98e4}  & \texttt{6.25e4} & \texttt{2.53e-5} & \texttt{1.58e0} & 99 & \,\ref{fig1} middle\\ 
 (\ref{sp1},\,a), $T_2$ & \texttt{2e-3}
 & \texttt{1.15e5}   & \texttt{1.21e5}  & \texttt{1.31e-5} & \texttt{1.58e0} & 258 &\,\ref{fig1} right \\ 
 (\ref{sp1},\,b), $T_3$ & \texttt{1e-3}
 & \texttt{1.07e0}  & \texttt{3.53e2}  & \texttt{2.82e-3} & \texttt{9.94e-1} & 246 & \,\ref{fig2} left \\
 (\ref{sp1},\,c), $T_3$ & \texttt{1e-3}
 & \texttt{5.78e1}  & \texttt{1.40e7} & \texttt{7.98e-8} & \texttt{1.12e0} & 27 & \,\ref{fig2} middle \\ 
 (\ref{sp1},\,c), $T_2$ & \texttt{1e-2}
 & \texttt{4.75e2}   & \texttt{1.15e8}  & \texttt{1.12e-8} & \texttt{1.29e0} & 171 & \,\ref{fig2} right  \\
\hline
\end{tabular}
\end{table}

\begin{remark}
As reported in~\cite[Fig.~2(b)]{SCB2025}, the example (\ref{sp1},\,a) is a problematic case for LOPCG
whose number of required steps exceeds $1000$ in a test where EPIC converges
with about $500$ steps. Here we consider similar cases and add more schemes
to the comparison. For the applied preconditioners, Table~\ref{tab1}
shows the associated \texttt{ichol} droptol values and quality parameters
$\kappa$, $\eta$ defined in \eqref{prec}, \eqref{thmpcge} for $\sigma=0$.
We additionally list the values $\beta_{\min}=\alpha_1\psi_2$
and $\beta_{\max}=\alpha_n\psi_n$ from the derivation \eqref{thmpcgf} of $\eta$
since $2\beta_{\min}$ and $2\beta_{\max}$ are close to
theoretical settings of the parameters $\mu$ and $L$ in EPIC and RAP.
\end{remark}

\begin{figure}[htbp]
\begin{center}
\includegraphics[trim={4.5cm 14.5cm 4.5cm 4.5cm}, clip, width=\textwidth]{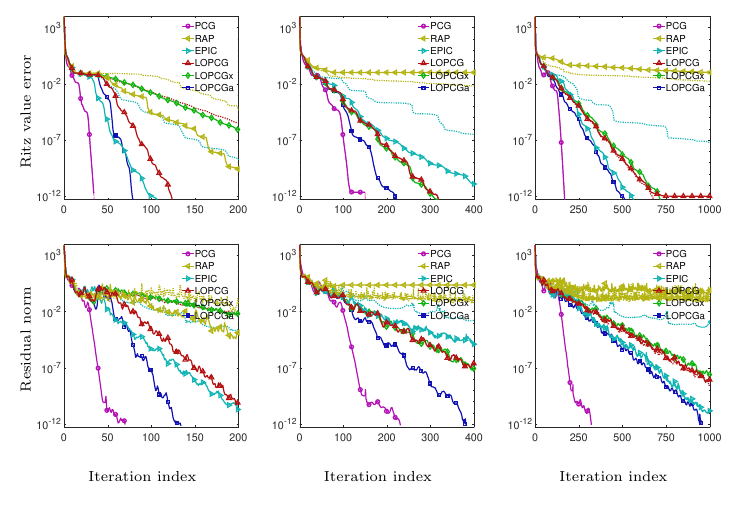}
\end{center}
\par\vskip -2ex
\caption{\small Numerical comparisons between CG-like eigensolvers
for the example (\ref{sp1},\,a); see Section~\ref{sec2c} for details.
The subplots include three cases where LOPCG implementations
are entirely/not/partially inferior to EPIC implementations. The augmented scheme LOPCGa
(Algorithm~\ref{alopcga}) is clearly more efficient.}
\label{fig1}
\end{figure}

We first implement LOPCG in the form of Algorithm~\ref{alopcg}
for various initial guesses with random number seeds among $\{0,\ldots,300\}$.
In this way we find a difficult case for LOPCG where it requires at least twice as many steps as EPIC. 
For completeness we also test two other versions of LOPCG. A version with full $M$-orthonormalization
following \cite{HL2006} shows no visible difference to Algorithm~\ref{alopcg},
whereas the code \cite[\texttt{lobpcg.m}]{KALO2007} does not suffer from
this troublesome case; see the left column in Figure~\ref{fig1}
for the convergence history. Therein the solid curves for LOPCG, EPIC and RAP
correspond to \texttt{lobpcg.m}, EPIC(6,\,6) and RAP(8,\,50) with empirical values
of the parameters $\mu$ and $L$ from~\cite{SCB2025,SC2025},
and the dotted curves stand for Algorithm~\ref{alopcg},
EPIC($2\beta_{\min},\,2\beta_{\max}$) and RAP($2\beta_{\min},\,2\beta_{\max}$).
Interestingly, although the empirical settings for EPIC and RAP
seem inconsistent with the theoretical construction, they lead to
better performance. For \texttt{lobpcg.m}, by looking into the details,
we note that the acceleration is enabled by switching to PSD in the $40$th step
where the condition number of the underlying Gram matrix exceeds a threshold.
In our additional test motivated by the unusual oscillation of the residual norm around that step,
a similar effective switch is made by observing the residual reduction.
Furthermore, the simple extension LOPCGx of Algorithm~\ref{alopcg}
(by adding $x^{(i-2)}$ to the trial subspace) only enables a slight acceleration.
In contrast, the adaptive extension LOPCGa (Algorithm~\ref{alopcga}) begins to accelerate
only a few steps later than EPIC and gets a much better convergence rate in the final phase.
The advantage of LOPCGa is more evident for the residual norm. We also illustrate the heuristic PCG
(Algorithm~\ref{apcg}) to indicate the space for improvement in further modifications.

The middle column in Figure~\ref{fig1} stands for a comparison
where the initial guess is selected such that the convergence speed of LOPCG is average among random tests.
In this case, the code \texttt{lobpcg.m} and Algorithm~\ref{alopcg}
have quite similar performance (so that the red dotted curve is covered)
and only slightly differ from LOPCGx. Moreover,
EPIC(6,\,6) is still competitive and superior to EPIC($2\beta_{\min},\,2\beta_{\max}$),
but slower than LOPCG in the final phase. The two implementations of RAP
behave like PSD, but RAP($2\beta_{\min},\,2\beta_{\max}$) becomes the better version.
The benefit of LOPCGa is substantial in the sense of closeness to the heuristic PCG.

The right column in Figure~\ref{fig1} is concerned with an even less accurate preconditioner and
the worst case of LOPCG among random tests. We note that the number of required steps
of EPIC(6,\,6) is comparable with that in \cite[Fig.~2(b)]{SCB2025}.
The two versions of LOPCG are not dramatically slower than EPIC(6,\,6)
and not accelerated by LOPCGx, whereas LOPCGa is again the closest one to the heuristic PCG.
The other schemes are far less competitive. A further fact is that
the expense per step of LOPCGa is slightly more than that of LOPCG and less than that of EPIC
since the steps are dominated by matrix-vector multiplications (with $A$, $M$ and $T$).
Thus LOPCGa also wins with respect to the total computational time.

\begin{figure}[htbp]
\begin{center}
\includegraphics[trim={4.5cm 14.5cm 4.5cm 4.5cm}, clip, width=\textwidth]{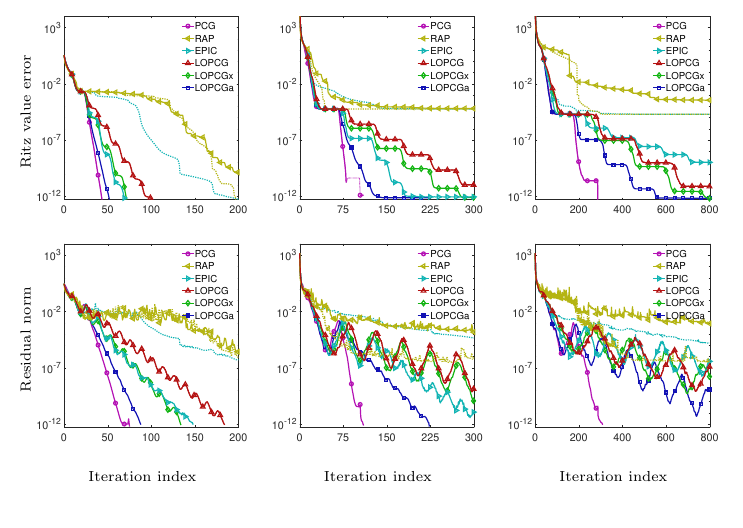}
\end{center}
\par\vskip -2ex
\caption{\small Numerical comparisons between CG-like eigensolvers
for the examples (\ref{sp1},\,b,\,c); see Section~\ref{sec2c} for details.
The subplots particularly indicate the cluster robustness of LOPCGa (Algorithm~\ref{alopcga})
in the sense that the residual norm reduction is straightforward or oscillates more weakly.}
\label{fig2}
\end{figure}

Figure~\ref{fig2} shows comparisons for the examples (\ref{sp1},\,b,\,c)
concerning clustered eigenvalues and the worst case of LOPCG among random tests.
It is remarkable that LOPCGx becomes more robust and mostly accelerates LOPCG.
Moreover, EPIC(6,\,6) is still the best one among the versions of EPIC and RAP or,
in more detail, EPIC mostly benefits from empirical settings, but RAP rarely.
Finally, the pleasing efficiency of LOPCGa motivates us to investigate
whether it can be further improved by more reliable augmentation strategies.

\section{Convergence theory of LOPCG}\label{sec3}

Despite the fairly early beginning of convergence analyses of CG-like eigensolvers
including a concept of LOPCG \cite[Section~1.5]{KNY1987}, direct estimates are only available
for certain schemes related to LOPCG, i.e., on the basis of comparative investigations.
The first-ever scheme is called an abstract two-stage method \cite[Section~3.2]{KNY1987}
and has with our settings the form
\begin{equation}\label{tsm}
x^{(i+1)}=f\big(T(A-\la^{(i)}M)\big)\,x^{(i)}
\end{equation}
where $f$ denotes an abstract matrix function. Although this scheme does not really cover
LOPCG due to the incompatible CG-term $p^{(i)}$, it corresponds to the weaker scheme PSD
(i.e. LOPCG without $p^{(i)}$) by specifying $f$ as a linear polynomial
so that indirect pessimistic estimates for LOPCG can be obtained. Moreover, if $f$ is specified
as a shifted Chebyshev polynomial, \eqref{tsm} can be interpreted as an acceleration
of PSD or a Lanczos-type scheme. This asymptotically matches multiple steps
of the generalized Davidson method in the final phase where $\la^{(i)}$ only varies slightly,
or of the heuristic PCG for sufficiently small $\la^{(i)}-\la_1$, leading to
indirect optimistic estimates for LOPCG. The framework of such estimates reads
\begin{equation}\label{fest}
\frac{\la^{(i+1)}-\la_1}{\la_2-\la^{(i+1)}}
\le\xi\,\frac{\la^{(i)}-\la_1}{\la_2-\la^{(i)}}
\end{equation}
under the assumption $\la^{(i)}<\la_2$. The convergence factor $\xi$
can be given in a cumbersome form involving $\la^{(i)}$
or in a concise form without $\la^{(i)}$. In particular, 
$\xi$ in the Chebyshev specification of \eqref{tsm}
is comparable with $\big(\mathcal{C}_i(\varphi)\big)^{-2}$ from
the estimate~\eqref{thmpcge} for the heuristic PCG.

Nevertheless, using optimistic estimates is sometimes problematic, especially
in the troublesome cases for LOPCG mentioned in Section~\ref{sec2c}.
In addition, deriving reasonable pessimistic estimates for LOPCG, e.g.,
via sharp estimates for PSD with Rayleigh-Ritz step sizes, was technically more challenging so that
the analysis in~\cite[Section~3.4]{KNY1987} actually arrived at a weaker
convergence factor which can be written as $\xi=1-\eta^{-1}$
by using the condition number $\eta$ defined in~\eqref{thmpcge}.
The next notable estimate for PSD is presented
in~\cite[Section~6]{OVT2006s} with $\xi$ depending on $\la^{(i)}$,
whereas the final form $\xi=\big((\eta-1)/(\eta+1)\big)^2$ expected in~\cite{KNY1987}
was fully proved in~\cite{NEY2012} after a number of investigations of PSD
with simple step sizes~\cite{KNN2003}. 
Here we reformulate the main result from~\cite{NEY2012} for further discussion.

\begin{lemma}\label{lmpsd}
Consider an arbitrary $x\in\CC^n{\setminus}\{0\}$ that is not an eigenvector of $(A,M)$
and let $\theta'$ be the minimum of the Rayleigh quotient \eqref{rq}
in $\spanl\{x,\,T(Ax-\theta\,Mx)\}$ for $\theta=\la(x)$. Then $\theta'<\theta$, and
\begin{equation}\label{lmpsde}
 \frac{\theta'-\la_j}{\la_{j+1}-\theta'}
 \le\xi_j\,\frac{\theta-\la_j}{\la_{j+1}-\theta}\quad\mbox{for}\quad
 \xi_j=\Big(\dfrac{\eta_j-1}{\eta_j+1}\Big)^2, \quad
 \eta_j=\kappa\,\frac{(\la_n-\la_j)(\la_{j+1}-\sigma)}{(\la_{j+1}-\la_j)(\la_n-\sigma)}
\end{equation}
with respect to the setting \eqref{prec} for preconditioners provided that $\la_j\le\theta<\la_{j+1}$.
\end{lemma}

The original formulation in~\cite{NEY2012} describes the preconditioner $T$
with a parameter corresponding to $(\kappa-1)/(\kappa+1)$ by \eqref{prec}.
The analysis therein is motivated by accelerating a perturbed shift-and-invert eigensolver
rather than comparing PSD with a gradient iteration
for solving the homogeneous linear system \eqref{lse}.

\begin{remark}
Applying Lemma~\ref{lmpsd} to PSD by considering $\theta=\la^{(i)}$
and $\theta'=\la^{(i+1)}$ indicates that the vector iterates
can at least reach an eigenvector associated with an interior eigenvalue $\la_j$
where a stagnation at $\la_j$ without eigenvector convergence is excluded by
the simple inequality $\theta'<\theta$. The convergence toward $\la_1$
cannot be ensured by Lemma~\ref{lmpsd} for arbitrary $T$ from~\eqref{prec}
without further assumptions due to the fact that
there exist examples fulfilling the equality in~\eqref{lmpsde}.
A technical reason is that the underlying analysis is eigenvalue-oriented
and related to a three-dimensional subspace concerning the slowest convergence.
An exception is the special case $T=A_{\sigma}^{-1}$ where
the convergence toward $\la_1$ can be shown on the basis of angle estimates
like~\cite[Theorem~1.2]{KNY1987}. Furthermore, Lemma~\ref{lmpsd}
provides intermediate bounds for LOPCG whose trial subspace is a superset
of that of PSD so that inequalities between eigenvalue approximations
are obtained by Courant-Fischer principles.
In contrast, angle estimates for PSD are not (directly) applicable to LOPCG.
\end{remark}

\begin{figure}[htbp]
\begin{center}
\includegraphics[trim={4.5cm 14.5cm 4.5cm 4.5cm}, clip, width=\textwidth]{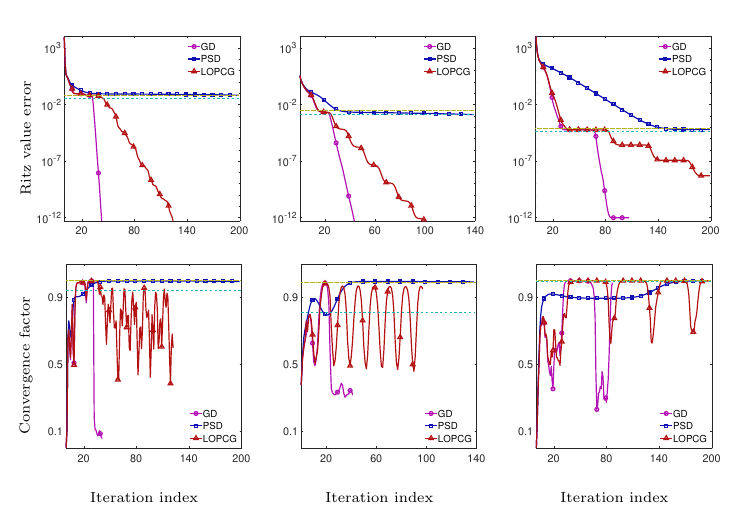}
\end{center}
\par\vskip -2ex
\caption{\small Convergence phenomena beyond the available estimates for
 CG-like eigensolvers (Remark~\ref{rmka}).
 Upper row: convergence history of $\la^{(i)}-\la_1$ in comparison to
 $\la_2-\la_1$ and $(\la_2-\la_1)/2$ (horizontal lines) which roughly
 mark the beginning of the final phase.
 Lower row: single-step convergence factor
 $(\la^{(i+1)}-\la_1)/(\la^{(i)}-\la_1)$ in comparison to
 the estimated values $\xi$ and $\psi^2$ (horizontal lines)
 from \eqref{lmpsde1} and \eqref{lopcge}.}
\label{fig3}
\end{figure}

\begin{remark}\label{rmka}
The convergence factor $\xi_j$ in~\eqref{lmpsde}
does not depend on any iterates so that the estimate can easily be applied to multiple steps.
In particular, \eqref{lmpsde} results in the following estimate concerning the final phase of PSD
with $\la^{(i)}<\la_2$:
\begin{equation}\label{lmpsde1}
 \frac{\la^{(i+m)}-\la_1}{\la_2-\la^{(i+m)}}
 \le\xi^m\,\frac{\la^{(i)}-\la_1}{\la_2-\la^{(i)}}\quad\mbox{for}\quad
 \xi=\Big(\dfrac{\eta-1}{\eta+1}\Big)^2, \quad
 \eta=\kappa\,\frac{(\la_n-\la_1)(\la_2-\sigma)}{(\la_2-\la_1)(\la_n-\sigma)}.
\end{equation}
With the same settings, the first-ever convergence estimate for PSD by Samokish \cite{SAM1958}
can be rewritten as
\begin{equation}\label{psde}
 \la^{(i+1)}-\la_1\le\xi(\la^{(i)}-\la_1)\big(1+\mathcal{O}(\sqrt{\la^{(i)}-\la_1}\,)\big),
\end{equation}
whereas \eqref{lmpsde1} with $m=1$ implies
\[\la^{(i+1)}-\la_1\le\xi(\la^{(i)}-\la_1)
 \Big(1+\frac{(\la^{(i)}-\la_1)-(\la^{(i+1)}-\la_1)}{(\la_2-\la_1)-(\la^{(i)}-\la_1)}\Big)
 =\xi(\la^{(i)}-\la_1)\big(1+\mathcal{O}(\la^{(i)}-\la_1)\big)\]
containing a slightly better asymptotic term. In contrast to this, the non-asymptotic estimate
$\la^{(i+1)}-\la_1\le\xi(\la^{(i)}-\la_1)$ does not hold in general and requires
a larger convergence factor instead of $\xi$; see \cite[Theorems~2.1 and 6.2]{OVT2006s}
with convergence factors depending on $\la^{(i)}$. It is still difficult to derive
a multi-step estimate like $\la^{(i+m)}-\la_1\le\tilde{\xi}^m(\la^{(i)}-\la_1)$
without using $\la^{(i)}$ to build an appropriate convergence factor~$\tilde{\xi}$.
Nevertheless, the asymptotic estimate \eqref{psde} can be upgraded
for the restarted generalized Davidson method as mentioned in~\cite[Eq.~(5.9)]{OVT2008a}.
This corresponds to
\begin{equation}\label{gde}
 \la^{(i+m)}-\la_1\le\big(\mathcal{C}_m(\varphi)\big)^{-2}(\la^{(i)}-\la_1)
 \big(1+\mathcal{O}(\sqrt{\la^{(i)}-\la_1}\,)\big)\quad\mbox{for}\quad
 \varphi=\frac{\eta+1}{\eta-1}
\end{equation}
with the settings for~\eqref{lmpsde1} and the Chebyshev polynomial $\mathcal{C}_m$.
It should be noted that \eqref{gde} also requires the assumption $\la^{(i)}<\la_2$
and cannot directly be derived on the basis of the similar estimate \eqref{thmpcge}
for the heuristic PCG. Furthermore, \eqref{lmpsde1} with
$\big(\mathcal{C}_m(\varphi)\big)^{-2}$ instead of $\xi^m$
corresponds to the conjectured estimate for LOPCG in~\cite{KNY2001}.
The multi-step convergence factor $\big(\mathcal{C}_m(\varphi)\big)^{-2}$
can also be transformed into
\begin{equation}\label{lopcge}
 \big(\mathcal{C}_m(\varphi)\big)^{-2}=\Big(\frac{2\psi^m}{1+\psi^{2m}}\Big)^2
 \quad\mbox{for}\quad\psi=\frac{\sqrt{\eta}-1}{\sqrt{\eta}+1}
\end{equation}
and then interpreted by the average single-step convergence factor $\psi^2$.
However, this does not mean that replacing $\xi$ by $\psi^2$ in \eqref{psde}
provides a suitable single-step estimate for LOPCG.
In Figure~\ref{fig3}, we compare LOPCG with PSD and
the non-restarted generalized Davidson method (GD)
concerning the first, fourth and fifth cases from Table~\ref{tab1}.
The Ritz value error $\la^{(i)}-\la_1$ is illustrated in the upper row
where the two horizontal lines denote $\la_2-\la_1$ and $(\la_2-\la_1)/2$.
The lower row shows the associated convergence factor
$(\la^{(i+1)}-\la_1)/(\la^{(i)}-\la_1)$ together with two horizontal lines
corresponding to $\xi$ and $\psi^2$ defined in \eqref{lmpsde1} and \eqref{lopcge}.
Obviously, each method is slowed down as soon as $\la^{(i)}$ reaches $\la_2$.
At this moment, the convergence factor is ``reset'' to $\xi$ even for GD
and causes a stagnation. This phenomenon can repeatedly occur during further
stagnations of LOPCG so that its convergence rate cannot simply be bounded by $\psi^2$.
On the other hand, $\psi^2$ does not provide an appropriate estimation
for the convergence acceleration of GD or LOPCG out of stagnations,
especially for less accurate preconditioners and clustered eigenvalues
where $\psi^2$ is close to $\xi$; cf.~the left and right columns in Figure~\ref{fig3}.
\end{remark}

\begin{remark}\label{rmkb}
The recent studies~\cite{SCB2025,SC2025} related to the convergence theory of LOPCG
present estimates for two alternative methods EPIC and RAP. We restate these results
with the settings for~\eqref{lmpsde1} for a brief discussion.
The estimate \cite[Eq.~(5.15)]{SCB2025} for EPIC (after correcting a typo) corresponds to
\begin{equation}\label{eqepic}
 \la^{(i+m)}-\la_1\le2\big(1-\sqrt{\tilde{\eta}^{-1}}\big)^m(\la^{(i)}-\la_1)
 \quad\mbox{for}\quad
 \tilde{\eta}=\eta+\mathcal{O}\big(\la(q)-\la_1\big)
\end{equation}
depending on an auxiliary vector $q$, and the estimate \cite[Eq.~(3)]{SC2025}
for RAP can be written as
\begin{equation}\label{eqrap}
 \la^{(i+m)}-\la_1\le2\Big(1-\frac12\sqrt{\tilde{\eta}^{-1}}\Big)^m(\la^{(i)}-\la_1)
 \quad\mbox{for}\quad
 \tilde{\eta}\approx\kappa\,\frac{\la_2-\sigma}{\la_2-\la_1}
\end{equation}
where the exact definition of $\tilde{\eta}$ requires an auxiliary scalar $\rho$
between $\la_1$ and $(\la_1+\la_2)/2$.
The estimates \eqref{eqepic} and \eqref{eqrap}
can easily be compared with \eqref{gde} by using the transformation \eqref{lopcge}
and neglecting the auxiliary terms. The convergence factor in \eqref{gde}
is clearly smaller than those in \eqref{eqepic} and \eqref{eqrap} and even than
the corresponding squared values. In comparison to \eqref{eqepic},
the added term $\frac12$ in \eqref{eqrap} is problematic in the sense
that the single-step convergence factor is at least $\frac12$,
leading to overestimations for well-separated eigenvalues and more precise preconditioners.
Moreover, \eqref{eqepic} and \eqref{eqrap} are weaker asymptotic estimates
since their auxiliary terms cannot be interpreted as monotonic decreasing
like $\mathcal{O}(\sqrt{\la^{(i)}-\la_1})$, e.g., the last auxiliary vector $q$
of EPIC is mostly chosen before the final phase due to \cite[Fig.~2]{SCB2025}.
\end{remark}

\begin{figure}[htbp]
\begin{center}
\includegraphics[trim={5cm 19cm 5cm 4.5cm}, clip, width=\textwidth]{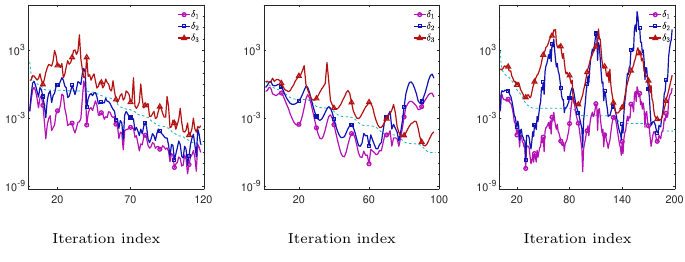}
\end{center}
\par\vskip -2ex
\caption{\small Limitation of asymptotic estimates for LOPCG (Remark~\ref{rmkc}).
 The asymptotic terms $\delta_1$, $\delta_2$, $\delta_3$
 for \eqref{lmlopcge}, \eqref{lmlopcge1}, \eqref{gde}
 can substantially deviate from the estimated size
 $\mathcal{O}(\sqrt{\la^{(i)}-\la_1}\,)$ for clustered eigenvalues
 (middle and right subplots).}
\label{cfc}
\end{figure}

\begin{remark}\label{rmkc}
The comparison in Remark~\ref{rmkb}, especially the exponent difference between the convergence factors,
indicates that LOPCG cannot thoroughly be investigated
by comparative approaches concerning more general problem settings.
We further note that some potential drawbacks of applying
the theory of non-quadratic CG minimization to CG-like eigensolvers
have been pointed out in \cite[Section~4]{OVT2008a}, followed by
the derivation of a direct two-step estimate for LOPCG on the basis of
asymptotically equivalent two-term CG-schemes. The explicit estimate essentially coincides with
the special form of \eqref{gde} for $m=2$ and can indeed be reconstructed by considering
that the decelerated form
\begin{equation}\label{gd2}
 x^{(i+1)} \ \ \xleftarrow{\mathrm{RRw}} \ \ x^{(i)}+\spanl\{Tr^{(i)},\,p^{(i)}\};
 \quad  x^{(i+2)} \ \ \xleftarrow{\mathrm{RRw}} \ \ x^{(i+1)}+\spanl\{Tr^{(i+1)}\};
\end{equation}
of LOPCG in its every second step corresponds to a two-step restarted version of the generalized Davidson method.
Furthermore, \cite[Theorem~4.1]{OVT2008a} provides a remarkable asymptotic estimate
for LOPCG: in the final phase with $\la^{(i)}<(\la_1+\la_2)/2$, it holds that
\begin{equation}\label{lmlopcge}
 (\la^{(i)}-\la^{(i+1)})^{-1}+(\la^{(i+1)}-\la^{(i+2)})^{-1}
 =(\la^{(i+1)}-\tilde{\la}^{(i+2)})^{-1}\big(1+\mathcal{O}(\sqrt{\la^{(i)}-\la_1}\,)\big)
\end{equation}
with the approximate eigenvalue $\tilde{\la}^{(i+2)}$ given by one PSD step applied to $x^{(i+1)}$.
The explicit estimate in \cite[Theorem~4.2]{OVT2008a} is based on
\begin{equation}\label{lmlopcge1}
 \frac{\la^{(i+1)}-\la_1}{\la^{(i)}-\la^{(i+1)}}+\frac{\la^{(i+2)}-\la_1}{\la^{(i+1)}-\la^{(i+2)}}
 =\frac{\tilde{\la}^{(i+2)}-\la_1}{\la^{(i+1)}-\tilde{\la}^{(i+2)}}+\mathcal{O}(\sqrt{\la^{(i)}-\la_1}\,)
\end{equation}
which is a transformation of \eqref{lmlopcge}, namely,
both sides of \eqref{lmlopcge} are multiplied by $\la^{(i+1)}-\la_1$ and then decreased by $1$.
However, this approach is somewhat problematic for clustered eigenvalues.
In Figure~\ref{cfc}, we consider again the first, fourth and fifth cases from Table~\ref{tab1}
together with the form $\delta=\mathcal{O}(\sqrt{\la^{(i)}-\la_1}\,)$ associated with the above estimates.
The corresponding $\delta_1$ for \eqref{lmlopcge} and $\delta_2$ for \eqref{lmlopcge1}
are compared with $\sqrt{\la^{(i)}-\la_1}$. In contrast to the left subplot with well-separated eigenvalues,
the other two subplots indicate that $\delta_1$ and $\delta_2$ can clearly exceed $\sqrt{\la^{(i)}-\la_1}$
after several stagnations caused by clustered eigenvalues.
In addition, we note that the approach from \cite{OVT2008a} could be extended by using
\[\delta_3=\big((\la^{(i)}-\la_1)^{-1/2}+(\la^{(i+2)}-\la_1)^{-1/2}
 -2(\tilde{\la}^{(i+2)}-\la_1)^{-1/2}\big)^{-1}\]
via three-term recurrences of Chebyshev polynomials to achieve the form \eqref{gde}.
Nevertheless, also $\delta_3$ is not comparable with $\sqrt{\la^{(i)}-\la_1}$
for clustered eigenvalues as seen in Figure~\ref{cfc}.
\end{remark}

Summarizing the above, the existing estimates
are meaningful for LOPCG concerning non-clustered eigenvalues,
yet the staircase-shaped convergence behavior of LOPCG in the case $\la_1\approx\la_2$
still cannot be fully explained due to $\big(\mathcal{C}_m(\varphi)\big)^{-2}\approx1$.
In particular, the acceleration between two delays needs to be investigated
by some other approaches in future research.

\section{Two-term PCG (TPCG) schemes}\label{sec4}

We turn attention to the two-term alternatives for LOPCG from~\cite{OVT2008,OVT2008a}
which can be formulated as direct modifications of Algorithm~\ref{alopcg}
as in Algorithm~\ref{apcgb}.

\begin{algorithm}
\caption{TPCG for computing $(\la_1,x_1)$ in \eqref{evp}}\label{apcgb}
 generate $x^{(0)}$; set $p^{(0)}=0$\;
\For{$i=0,1,\ldots$ {\rm until convergence}}
{$\la^{(i)}=\la(x^{(i)})$; \quad $r^{(i)}=Ax^{(i)}-\la^{(i)}Mx^{(i)}$\;
 select $\iota^{(i)},\tau^{(i)}$; \quad
 $p^{(i+1)}=\iota^{(i)}\,Tr^{(i)} + \tau^{(i)} p^{(i)}$\;
 $x^{(i+1)} \ \ \xleftarrow{\mathrm{RRw}} \ \ x^{(i)}+\spanl\{p^{(i+1)}\}$ as in \eqref{rrw}\;}
\end{algorithm}

Therein the essential change is that the update of the search direction $p^{(i)}$
does not depend on the Rayleigh-Ritz procedure. Moreover, according to
our numerical tests, the weighted form $x^{(i+1)}=x^{(i)}+\cdots$
is particularly appropriate for TPCG without frequent normalizations
since it avoids overflows of $\|x^{(i+1)}\|_M$. The variants of Algorithm~\ref{apcgb}
only differ in the scalar parameters $\iota^{(i)}$ and $\tau^{(i)}$. For instance,
with an auxiliary $\gamma^{(-1)}=1$, one can define
\begin{equation}\label{apcgbg}
\begin{array}{cl}
\mbox{(a)} & \iota^{(i)}=2\,\|x^{(i)}\|_M^{-2} \quad\mbox{and}\quad
 \tau^{(i)}=\gamma^{(i)}/\gamma^{(i-1)} \quad\mbox{for}\quad
 \gamma^{(i)}=(\iota^{(i)})^2\,\|r^{(i)}\|_T^2, \\[1ex]
\mbox{(b)} &  \mbox{alternative for (a)} \quad\mbox{with}\quad
 \tau^{(i)}=\big(\gamma^{(i)}-(\iota^{(i)}r^{(i)})^*T
 (\iota^{(i-1)}r^{(i-1)})\big)/\gamma^{(i-1)}, \\[1ex]
\mbox{(c)} & \iota^{(i)}=1, \quad  \tau^{(0)}=1 \quad\mbox{and}\quad
 \tau^{(i)}=-\big(p^{(i)}{}^*BTr^{(i)}\big)/\big(p^{(i)}{}^*Bp^{(i)}\big)
 \ \ \mbox{for} \ \ i>0 \\[1ex]
 & \mbox{with} \ \ B=Q_{\alpha}^*(A-\beta M)Q_{\alpha}, 
 \ \ Q_{\alpha}=I - \alpha \|x^{(i)}\|_M^{-2}\,x^{(i)}x^{(i)}{}^*M
 \ \ \mbox{and} \ \ \beta\le\la^{(i)}.
\end{array}
\end{equation}
In particular, (\ref{apcgbg},\,a) can be regarded as
a modification of the heuristic PCG Algorithm~\ref{apcg}
where the gradient of the quadratic function $f(x)=\frac12(x^*A_{\la_1}x)$
is replaced by that of the Rayleigh quotient \eqref{rq}, i.e., $\iota^{(i)}r^{(i)}$.
Furthermore, (\ref{apcgbg},\,a) and (\ref{apcgbg},\,b)
correspond to the Bradbury-Fletcher and Polak-Ribi\'ere schemes
mentioned in~\cite[Section~3]{OVT2008a}, whereas (\ref{apcgbg},\,c)
matches the Daniel (with $\alpha=2$, $\beta=\la^{(i)}$),
Perdon-Gambolati (with $\alpha=0$, $\beta<\la_1$)
and Jacobi (with $\alpha=1$, $\beta\le\la^{(i)}$) schemes therein, and directly
defines the conjugation between two consecutive search directions by a $B$-orthogonality.
The matrix $B$ is a scaled Hessian of the Rayleigh quotient \eqref{rq} in the Daniel scheme
and an approximate Hessian of the quadratic function $f(x)=\frac12(x^*A_{\la_1}x)$
in the Perdon-Gambolati scheme.

The Jacobi scheme is less related to the Hessian of an objective function,
but rather to the Jacobi orthogonal complement correction (JOCC) equation \cite{SV1996}.
We briefly introduce the JOCC equation with the settings for Algorithm~\ref{apcgb}:
by neglecting the trivial case that $x^{(i)}$ is collinear with
or $M$-orthogonal to a target eigenvector $x_1$,
any nonzero vector $y\in\spanl\{x_1,x^{(i)}\}\cap\spanl\{x^{(i)}\}^{\perp_M}$
(where the superscript ${}^{\perp_M}$ denotes $M$-orthogonal complements)
can be scaled as $\tilde{x}^{(i)}=-x^{(i)}+\tau x_1$ so that
\[(A-\la_1M)\tilde{x}^{(i)}=-(A-\la_1M)x^{(i)}+\tau(A-\la_1M)x_1
 =-r^{(i)}-(\la^{(i)}-\la_1)Mx^{(i)}.\]
By using the projection matrix of $\spanl\{x^{(i)}\}^{\perp_M}$,
i.e., $Q^{(i)}=Q_1$ from (\ref{apcgbg},\,c), and the property
$\tilde{x}^{(i)}=Q^{(i)}\tilde{x}^{(i)}$, one gets the JOCC equation
\begin{equation}\label{jocc}
 J_{\la_1}^{(i)}\tilde{x}^{(i)}=-r^{(i)}
 \quad\mbox{for}\quad J_{\la_1}^{(i)}=Q^{(i)}{}^*(A-\la_1M)Q^{(i)}, \quad
 Q^{(i)}=I - \|x^{(i)}\|_M^{-2}\,x^{(i)}x^{(i)}{}^*M.
\end{equation}
The JOCC matrix $J_{\la_1}^{(i)}$ in~\eqref{jocc} is evidently Hermitian positive semidefinite
so that one can heuristically solve \eqref{jocc} by a PCG iteration to determine
the eigenvector $x_1=\tau^{-1}(x^{(i)}+\tilde{x}^{(i)})$. In practice,
the unknown eigenvalue $\la_1$ is replaced by a guess $\beta\le\la^{(i)}$
so that the search directions are conjugate with respect to
$J_{\beta}^{(i)}=Q^{(i)}{}^*(A-\beta M)Q^{(i)}$,
i.e., a special form of $B$ in  (\ref{apcgbg},\,c). Solving this practical version of \eqref{jocc}
is associated with inner steps of the Jacobi-Davidson CG  (JDCG) scheme \cite{NY2002}.
The Jacobi scheme from \cite{OVT2008a} can roughly be regarded as an alternative
for JDCG with degenerate inner steps where each $J_{\beta}^{(i)}$ is used
only once and acts as a connector of two outer steps.

Moreover, the notable analysis in~\cite[Section~3]{OVT2008a} shows that
the TPCG schemes except the Perdon-Gambolati scheme
are asymptotically equivalent to LOPCG in the final phase;
see also the numerical comparison between their block versions in~\cite{OVT2008}.
Thus it is possible to reduce the total computational time of LOPCG
by using TPCG. We prefer the Jacobi scheme since its construction
is closely associated with LOPCG so that it is reasonable to extend it
as an alternative for LOPCGa using 3D instead of 4D trial subspaces.
However, we note that the Jacobi scheme is implemented in~\cite{OVT2008}
with another approach for updating the search direction.
For Algorithm~\ref{apcgb}, this means that
$p^{(i+1)}$ is determined by a linear combination of $Tr^{(i)}$
and a non-target Ritz vector from the previous step instead of $p^{(i)}$.
This difference motivates us to formulate a compact derivation and investigate
whether a direct update using $p^{(i)}$ impedes proper performance.

The following derivation of the Jacobi scheme is a heuristic reformulation
of LOPCG based on \cite[Section~3.2]{OVT2008a}.
In Algorithm~\ref{alopcg}, the first two trial subspaces can be represented by
\[\mathcal{V}=\spanl\{x^{(0)},x^{(1)}\} \quad\mbox{and}\quad
 \mathcal{W}=\spanl\{x^{(0)},x^{(1)},Tr^{(1)}\}=\spanl\{x^{(0)},x^{(1)},x^{(2)}\}.\]
By using the projection matrix $Q^{(1)}$
defined in \eqref{jocc} and the search direction
$p^{(1)}=x^{(1)}-x^{(0)}\in\mathcal{V}$, we consider
$q^{(1)}=Q^{(1)}p^{(1)}\in\mathcal{V}\subset\mathcal{W}$.
Since $(\la^{(i)},x^{(i)})$ is a Ritz pair by construction, the fact that the Ritz vector residual
is orthogonal to the corresponding subspace, implies
\[0=q^{(1)}{}^*(A-\la^{(1)}M)x^{(1)}
 =q^{(1)}{}^*(A-\la^{(2)}M)x^{(1)} - (\la^{(1)}-\la^{(2)})q^{(1)}{}^*Mx^{(1)}
 =q^{(1)}{}^*(A-\la^{(2)}M)x^{(1)}\]
and $0=q^{(1)}{}^*(A-\la^{(2)}M)x^{(2)}$. Thus $(A-\la^{(2)}M)q^{(1)}$
is orthogonal to any linear combination of $x^{(1)}$ and $x^{(2)}$,
e.g., $Q^{(1)}x^{(2)}$. In other words, we get
\[0=(Q^{(1)}x^{(2)})^*(A-\la^{(2)}M)q^{(1)}
 =x^{(2)}{}^*y \quad\mbox{for}\quad
 y=Q^{(1)}{}^*(A-\la^{(2)}M)Q^{(1)}p^{(1)}\]
so that $x^{(2)}$ belongs to
$\mathcal{U}=\mathcal{W}\cap\spanl\{y\}^{\perp}$. Apart from trivial cases,
$\mathcal{U}$ has dimension $2$ and contains $x^{(1)}$ due to $Q^{(1)}x^{(1)}=0$.
A further vector in $\mathcal{U}$ can be obtained by orthogonalizing
$Tr^{(1)}\in\mathcal{W}$ against $y$ and acts as the next search direction $p^{(2)}$, namely,
\begin{equation}\label{apcgc}
 p^{(2)}=Tr^{(1)}-\frac{y^*Tr^{(1)}}{y^*p^{(1)}}\,p^{(1)}
 =Tr^{(1)}-\frac{p^{(1)}{}^*Q^{(1)}{}^*(A-\la^{(2)}M)Q^{(1)}Tr^{(1)}}
 {p^{(1)}{}^*Q^{(1)}{}^*(A-\la^{(2)}M)Q^{(1)}p^{(1)}}\,p^{(1)}.
\end{equation}
Then applying the Rayleigh-Ritz procedure to $\spanl\{x^{(1)},p^{(2)}\}$ yields $x^{(2)}$,
and the reduced trial subspace $\spanl\{x^{(1)},p^{(2)}\}=\spanl\{x^{(1)},x^{(2)}\}$
can be used to reformulate the next step of LOPCG in the same way. This successive reformulation
leads to Algorithm~\ref{apcgb} where the counterpart of \eqref{apcgc}
for the $(i{+}1)$th step turns into the update formula (\ref{apcgbg},\,c)
with $\alpha=1$ and a shift $\beta\le\la^{(i)}$ instead of $\la^{(i+1)}$.

Furthermore, since the reduced trial subspaces are two-dimensional,
$q^{(1)}=Q^{(1)}p^{(1)}$ and its counterparts $q^{(i)}=Q^{(i)}p^{(i)}$ in further steps
are actually non-target Ritz vectors as used in the implementation in~\cite{OVT2008}.
The purpose of using these Ritz vectors is likely to bypass the orthogonalization
in $Q^{(i)}p^{(i)}$. However, we note that $Q^{(i)}p^{(i)}$ can be computed
in an equivalent form with ready-to-use scalars so that a direct update using $p^{(i)}$
is still efficient. More precisely, for
$Q^{(i)}p^{(i)}=p^{(i)}- \|x^{(i)}\|_M^{-2}\,x^{(i)}x^{(i)}{}^*Mp^{(i)}$,
one can evaluate the dot product
\[x^{(i)}{}^*Mp^{(i)}=(x^{(i-1)}+\delta p^{(i)})^*Mp^{(i)}
 =x^{(i-1)}{}^*Mp^{(i)}+\delta p^{(i)}{}^*Mp^{(i)}\]
by using components of the Gram matrix $U^*MU$ for $U=[x^{(i-1)},\,p^{(i)}]$
and the step size $\delta$ from the previous step.
Thus $Q^{(i)}p^{(i)}$ is obtainable by a simple vector update.

As another important factor for practical implementation,
the shift $\beta\le\la^{(i)}$ needs to be specified.
We note that the choice $\beta=\la^{(i)}$ in~\cite{OVT2008}
is not always close to the heuristic value $\la^{(i+1)}$, especially in the starting phase.
In this sense, setting $\beta<\la^{(i)}$ is more meaningful, e.g.,
\begin{equation}\label{apcgs}
 \beta=\max\big\{(\sigma+\la^{(i)})/2, \,\ 2\la^{(i)}-\la^{(i-1)}\big\}
\end{equation}
with an estimated lower bound $\sigma<\la_1$ concerning the preconditioner \eqref{prec}.
Interestingly, this shift is also appropriate for a slightly modified search direction update
using $Q^{(i-1)}p^{(i)}$ (inspired by similar algorithmic simplifications like Gauss-Seidel).
This combination does not deteriorate, but even slightly improves the overall performance.

\begin{algorithm}
\caption{Two specified TPCG schemes for computing $(\la_1,x_1)$ in \eqref{evp}}\label{atpcg}
 generate $x^{(0)}$\;
\For{$i=0,1,\ldots$ {\rm until convergence}}
{$\la^{(i)}=\la(x^{(i)})$; \quad $r^{(i)}=Ax^{(i)}-\la^{(i)}Mx^{(i)}$\;
 \eIf{$i>0$}
 {$v=Q^{(i)}p^{(i)}$ \ or \ $Q^{(i-1)}p^{(i)}$\;
 $w=Av-\beta Mv$ \ with $\beta$ defined in \eqref{apcgs}\;
 $\tau^{(i)}=-(w^*Tr^{(i)})/(w^*v)$\;
 $p^{(i+1)}=Tr^{(i)}+\tau^{(i)}p^{(i)}$ \ or \ $Tr^{(i)}+\tau^{(i)}v$\;}
 {$p^{(i+1)}=Tr^{(i)}$\;}
 $x^{(i+1)} \ \ \xleftarrow{\mathrm{RRw}} \ \ x^{(i)}+\spanl\{p^{(i+1)}\}$
 \ or \ $x^{(i)}+\spanl\{Q^{(i)}p^{(i+1)}\}$ as in \eqref{rrw}\;}
\end{algorithm}

The above discussion leads to two specified TPCG schemes in Algorithm~\ref{atpcg}
whose differences can be seen in Lines 5, 8 and 11. For the first (and standard) version,
one needs to note that the term $w^*Tr^{(i)}$
in Line 7 coincides with $w^*Q^{(i)}Tr^{(i)}$ and thus
does not contradicts the description in (\ref{apcgbg},\,c).
This coincidence directly follows from
\[w^*x^{(i)}=p^{(i)}{}^*Q^{(i)}{}^*(A-\beta M)x^{(i)}
 =p^{(i)}{}^*Q^{(i)}{}^*(A-\la^{(i)}M)x^{(i)}=p^{(i)}{}^*r^{(i)}=0\]
where the last equality uses the fact that $r^{(i)}$
is a Ritz vector residual of the previous trial subspace containing $p^{(i)}$.
Concerning the performance of Algorithm~\ref{atpcg}, the Rayleigh-Ritz procedure
applied to the reduced trial subspaces is merely an elementary calculation
so that the expense per step is only slightly more than that of the heuristic PCG
and less than that of LOPCG. The concrete difference depends on how strongly
the three indispensable explicit matrix-vector multiplications (with $A$, $M$ and $T$)
dominate the step (or equivalently, on the nonzero component density of these matrices).
For instance, if applying the second version of Algorithm~\ref{atpcg} to the cases
from Table~\ref{tab1} and using the expense of the heuristic PCG as the unit of measure,
then the expense of LOPCG or Algorithm~\ref{atpcg} varies between $1.18$ and $2.31$
or between $1.06$ and $1.39$. Therein LOPCG is tested
with \cite[\texttt{lobpcg.m}]{KALO2007}, yet a slight acceleration can be achieved
by skipping certain lines which are not necessary for single-vector iterations. 

\begin{figure}[htbp]
\begin{center}
\includegraphics[trim={5cm 19cm 5cm 4.5cm}, clip, width=\textwidth]{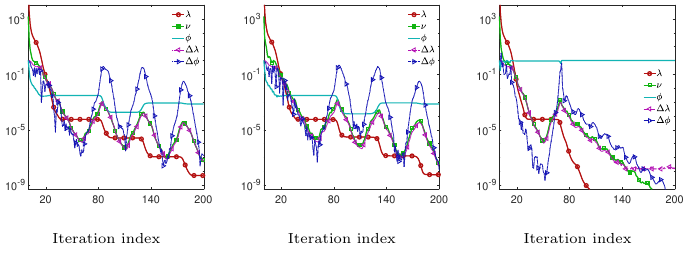}
\end{center}
\par\vskip -2ex
\caption{\small Comparing LOPCG and TPCG with LOPCGa
in terms of \eqref{spcgm} to suggest timely augmentation.
Changes of the residual norm $\nu$ can provide good timing.}
\label{spcg1}
\end{figure}

Algorithm~\ref{atpcg} can easily be upgraded to a scheme that is more efficient
than the winner LOPCGa in the numerical comparisons in Section~\ref{sec2c}.
The trial subspaces are also augmented by an auxiliary vector $a$,
yet the simple strategy $|\cos\angle_M(a,\,x^{(i)})|<0.7$ for updating $a$ in LOPCGa
is refined by observing residual norms motivated by the following test.
Therein the second version of Algorithm~\ref{atpcg} is denoted by TPCG
and compared with LOPCG and LOPCGa
concerning the fifth case from Table~\ref{tab1} (containing clustered eigenvalues).
As seen in Figure~\ref{spcg1}, their convergence history is illustrated in terms of
\begin{equation}\label{spcgm}
\begin{split}
 & \la: \ \la^{(i)}-\la_1,\quad \nu: \ \|r^{(i)}\|_2/\|x^{(i)}\|_M, \quad
 \phi: \ \phi^{(i)}=|\cos\angle_M(a,\,x^{(i)})| \\[1ex]
 & \hspace{1cm} \Delta\lambda: \ |\la^{(i)}/\la^{(i-1)}-1|^{1/2}, \quad
 \Delta\phi: \ |\phi^{(i)}/\phi^{(i-1)}-1|.
\end{split}
\end{equation}

This comparison clearly reflects the asymptotic equivalence
between LOPCG (left subplot) and TPCG (middle subplot).
Moreover, since the previous test in Figure~\ref{fig2} shows that the starting phase of LOPCG
cannot be significantly accelerated, we only need to observe the data as of the first delay.
Guided by the strategy $\phi^{(i)}=|\cos\angle_M(a,\,x^{(i)})|<0.7$ for LOPCGa (right subplot),
we note that substantial changes of $\phi^{(i)}$
for $a=x^{(0)}$ are almost simultaneous with peaks of $\nu$.
This relation is more obvious by observing $\Delta\phi$.
Therefore a more reasonable augmentation strategy is to detect the residual peaks.
Another remarkable phenomenon is that $\Delta\lambda$ almost matches with $\nu$,
and is probably useful for a novel approach in future investigation of LOPCG and TPCG.

\begin{figure}[htbp]
\begin{center}
\includegraphics[trim={5cm 19cm 5cm 4.5cm}, clip, width=\textwidth]{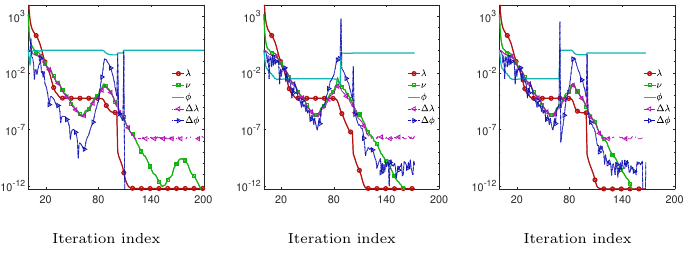}
\end{center}
\par\vskip -2ex
\caption{\small Comparing augmented versions of TPCG using
the angle-based criterion $\phi^{(i)}=|\cos\angle_M(a,\,x^{(i)})|<\tau$
or a residual-based criterion. Left: with $\tau=0.99$. Middle: with $\tau=0.001$.
Right: update when the residual norm $\nu$ decreases again after a substantial increase.}
\label{spcg2}
\end{figure}

We further note that the simple strategy $\phi^{(i)}<\tau$ with $\tau=0.7$ is not suitable
for the augmentation of TPCG and even causes a deceleration. The left and middle subplots
in Figure~\ref{spcg2} illustrate the acceptable performance by using two other thresholds:
$\tau=0.99$ and $\tau=0.001$. The latter is more effective than the former
in the sense of fewer updates in the starting phase and no further delays in the final phase.
Moreover, $a$ is initialized by $x^{(0)}$ and updated by the iterate corresponding to
the minimum of the (already computed) residual norm $\nu$. This mostly improves
the trivial update setting $a=x^{(i)}$ as concluded from our additional tests.
Next, the right subplot presents a convenient residual-based strategy which allows
reducing dependence on the threshold choice. In our attempts behind that,
three time points are tested to update the augmentation:
(i) directly after a substantial increase (valley) of the residual norm $\nu$;
(ii) when $\nu$ decreases again; (iii) when $\nu$ becomes smaller than its valley value.
The approaches (i) and (ii) well avoid the delays of the non-augmented TPCG or LOPCG.
However, (i) slightly slows down the uphill steps after the valley of $\nu$
so that a better timing is given by (ii). The approach (iii) is less appropriate due to
a further (but only one) delay after the update. Thus the approach (ii) is demonstrated
in Figure~\ref{spcg2} (right) where the hurdle-free downhill path of $\nu$ is comparable with
that by using the threshold $\tau=0.001$ in Figure~\ref{spcg2} (middle). In addition, $\nu$
falls below \texttt{1e-12} before the $160$th step which means a significant improvement
in comparison to the $225$th step of LOPCGa observed in Figure~\ref{fig2} (middle). 

\begin{figure}[htbp]
\begin{center}
\includegraphics[trim={5cm 15cm 5cm 4.5cm}, clip, width=\textwidth]{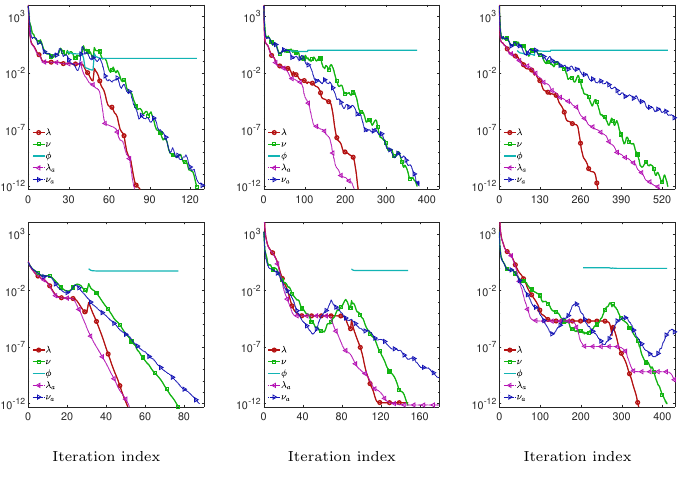}
\end{center}
\par\vskip -2ex
\caption{\small Comparing TPCGa with LOPCGa in terms of \eqref{spcgm} where $\la_a$
and $\nu_a$ stand for $\la$ (Ritz value error) and $\nu$ (residual norm) applied to LOPCGa.
The upper/lower row can be combined with the comparison in Figure \ref{fig1}/\ref{fig2}.
The benefit of TPCGa is significant for less accurate preconditioners and clustered eigenvalues.}
\label{spcg3}
\end{figure}

Summarizing the above, we suggest to modify Algorithm~\ref{atpcg}
by defining $\nu^{(i)}=\|r^{(i)}\|_2/\|x^{(i)}\|_M$ and adding to
\[\begin{array}{cl}
\mbox{Line 3 :} & \mbox{store the iterate $\check{x}$
 corresponding to $\nu_{\mbox{\tiny min}}=\min_{j \le i}\nu^{(j)}$}, \\[1ex]
\mbox{Line 11 :} & \mbox{augment the trial subspace by $\check{x}$
 if $\nu^{(i)}$ reaches a peak}
\end{array}\]
where the peak can be detected by getting \texttt{flag=2} from
\begin{equation}\label{atpcgt}
\begin{array}{ll}
\mbox{if \ \texttt{flag==0} \ and \ $\nu^{(i)}>1.5\,\nu_{\mbox{\tiny min}}$,}
 & \mbox{set \ \texttt{flag=1}}, \\[1ex]
\mbox{if \ \texttt{flag==1} \ and \ $\nu^{(i)}<\nu^{(i-1)}$,}
 & \mbox{set \ \texttt{flag=2}}.
\end{array}
\end{equation}
Since the augmented trial subspace has dimension $3$,
the expense per step of the resulting ``TPCGa'' scheme is comparable with that of LOPCG so that
the reduction of the total computational time is directly reflected by the number of required steps.
Moreover, we note that \eqref{atpcgt} is particularly appropriate for computing clustered eigenvalues.
In other cases such as the examples in Figure~\ref{fig1},
the residual norm could oscillate very frequently so that
a perfect timing cannot easily be predicted. Then it is more effective
to check a multi(e.g. $3$)-step decrease of $\nu$ instead of
the single-step decrease $\nu^{(i)}<\nu^{(i-1)}$.
Last but not least, the residual-based criterion allows us to drop the augmentation
in the starting phase so that the total expense can be further reduced. 
Following Figures \ref{fig1} and \ref{fig2}, we compare this final version of TPCGa
with LOPCGa in Figure~\ref{spcg3}. Therein $\la_a$ and $\nu_a$ denote the terms
$\la$ and $\nu$ from \eqref{spcgm} applied to LOPCGa.

\section*{Conclusion}

Competitive eigensolvers can be constructed by mimicking or modifying
the traditional CG iteration for solving linear systems. However, the optimality
of the original search directions is not preserved, but occasionally heavily disturbed
as observed in the staircase-shaped convergence history of LOPCG
and its recent alternatives while computing clustered eigenvalues.
Our single-vector scheme TPCGa with augmented two-term recurrences overcomes
this drawback by detecting peaks of residual norms and punctually renewing the augmentation.
The low dimension (two or three) of the trial subspace in TPCGa ensures that
its expense per step is even less than that of LOPCG, leading to
notably better performance also with respect to the total computational time
(not restricted to the case of clustered eigenvalues). This approach also induces
future topics concerning blockwise implementations and convergence analyses.

\appendix

\section{Typical generalized eigenvalue problems}

The introduction of LOPCG in \cite{KNY2001} is concerned with
computing the largest eigenvalue of
\begin{equation}\label{evpc}
 Mx=\mu Ax, \quad M,A\in\CC^{n \times n}\ \mbox{Hermitian}, \quad
 A\ \mbox{positive definite}.
\end{equation}
This covers two cases of the eigenvalue problem
\begin{equation}\label{evpd}
 Lu=\la Su, \quad L,S\in\CC^{n \times n}\ \mbox{Hermitian}, \quad
 S\ \mbox{positive definite}
\end{equation}
arising from a discretization of a self-adjoint elliptic partial differential operator.

The first case coincides with the case introduced
in Section~\ref{sec1} (by setting $L=A$, $S=M$),
i.e., the smallest eigenvalue $\la_1$ of the matrix pair $(L,S)$ is of interest.
For matching the form \eqref{evpc}, one uses a shift $\sigma<\la_1$
so that \eqref{evpd} can be reformulated as $\tilde{L}u=(\la-\sigma)Su$
with the shifted matrix $\tilde{L}=L-\sigma S$ which is evidently positive definite.
Then \eqref{evpc} corresponds to \eqref{evpd} 
by setting $M=S$, $A=\tilde{L}$ and $\mu=(\la-\sigma)^{-1}$.
Moreover, in a practical treatment with preconditioned eigensolvers,
linear systems of the form $\tilde{L}v=w$ are to be solved approximately.

The second case deals with a shift $\sigma$ close to but not equal to
an interior eigenvalue, i.e., $\tilde{L}$
is indefinite and invertible. Then \eqref{evpd} can be reformulated as
$(\tilde{L}S^{-1}\tilde{L})u=(\la-\sigma)\tilde{L}u$
and matches \eqref{evpc} with
$M=\pm\tilde{L}$, $A=\tilde{L}S^{-1}\tilde{L}$
and $\mu=\pm(\la-\sigma)^{-1}$. Therein the symbol $\pm$ means
determining eigenvalues either larger or smaller than $\sigma$.
An alternative representation uses
$M=S$, $A=\tilde{L}S^{-1}\tilde{L}$
and $\mu=(\la-\sigma)^{-2}$. Applying preconditioned eigensolvers
to this case requires approximately solving linear systems
of the form $(\tilde{L}S^{-1}\tilde{L})v=w$
which corresponds to two linear systems for the matrix $\tilde{L}$
in the implementation.

A further special form of \eqref{evpc} concerns Hermitian definite
matrix pencils \cite{KPS2014}
\begin{equation}\label{evpdp}
 Lu=\la Su, \quad L,S\in\CC^{n \times n}\ \mbox{Hermitian}, \quad
 \tilde{L}=L-\sigma S\ \mbox{definite for a certain}\ \sigma\in\RR.
\end{equation}
If $\tilde{L}$ is positive definite, the transformation into \eqref{evpc}
is the same as in the first case for \eqref{evpd}. For negative definite $\tilde{L}$,
one sets $M=S$, $A=-\tilde{L}$ and $\mu=-(\la-\sigma)^{-1}$ instead.
Moreover, changing the signs of $M$ and $\mu$ (in both cases)
allows determining eigenvalues on both sides of $\sigma$.
A practical eigenvalue arrangement with respect to the inertia of
the possibly indefinite and singular matrix $S$ is suggested in \cite[Theorem~2.1]{KPS2014}.

One can also use \eqref{evpc} for describing
the linear response eigenvalue problem~\cite{BL2012}.
Therein the smallest positive eigenvalues of the block anti-diagonal matrix
\begin{equation}\label{evplr}
 \begin{bmatrix}O & \tilde{L} \\[1ex] \tilde{S} & O\end{bmatrix},
 \quad \tilde{L},\tilde{S}\in\CC^{\tilde{n} \times \tilde{n}}\
 \mbox{Hermitian positive semidefinite, one of them definite}
\end{equation}
are to be computed. A simple reformulation uses the fact that the target eigenvalues
are square roots of the smallest nonzero eigenvalues of the matrix pair
$(\tilde{S},\tilde{L}^{-1})$ or $(\tilde{L},\tilde{S}^{-1})$
provided that $\tilde{L}$ or $\tilde{S}$ is definite.
This is evidently included in the first case for \eqref{evpd}.
If both of $\tilde{L}$ and $\tilde{S}$ are definite,
a reformulation suggested in \cite{KPS2014} reads
\begin{equation}\label{evplrp}
 Lu=\la Su \quad\mbox{with}\quad
 L=\begin{bmatrix}\tilde{L} & O \\[1ex] O & \tilde{S}\end{bmatrix}
 \quad\mbox{and}\quad
 S=\begin{bmatrix}O & \tilde{I} \\[1ex] \tilde{I} & O\end{bmatrix}
\end{equation}
where $\tilde{I}$ denotes the $\tilde{n}{\times}\tilde{n}$ identity matrix.
This matches \eqref{evpdp} with $\sigma=0$.
Transformations of \eqref{evpd} and \eqref{evpdp} into \eqref{evpc}
enable constructing suitable eigensolvers based on those for \eqref{evpc}.

\end{document}